\documentclass{amsart}

\usepackage{amsmath, amsthm}
\usepackage{amssymb}
\usepackage{subfiles}

\usepackage{cite}

\usepackage{tikz}
\usetikzlibrary{positioning}

\theoremstyle{definition}

\newtheorem{thm}{Theorem}[section]
\newtheorem*{thm*}{Theorem}
\newtheorem{lemma}[thm]{Lemma}
\newtheorem{lem}[thm]{Lemma}
\newtheorem{defn}[thm]{Definition}
\newtheorem{claim}[thm]{Claim}
\newtheorem{prop}[thm]{Proposition}
\newtheorem{cor}[thm]{Corollary}
\newtheorem{const}[thm]{Construction}
\newtheorem{remark}[thm]{Remark}
\newtheorem{fact}[thm]{Fact}

\newtheorem{ex}[thm]{Example}

\renewcommand{\subset}{\subseteq}
\renewcommand{\supset}{\supseteq}

\newcommand\Erel{\mathrel{E}}
\newcommand\Frel{\mathrel{F}}

\newcommand{\rest}{\restriction}

\newcommand\cce{\mathrm{CC}[E^\omega]}
\newcommand\ccf{\mathrm{CC}[F^\omega]}

\newcommand{\action}{\curvearrowright}

\newcommand\lto{\longrightarrow}
\newcommand\force{\Vdash}
\newcommand\R{\mathbb{R}}
\newcommand\N{\mathbb{N}}

\newcommand\Z{\mathbb{Z}}
\newcommand\F{\mathbb{F}_2}
\newcommand\es{\emptyset}

\newcommand\im{\mathrm{Im}}

\DeclareMathOperator{\dom}{dom}
\DeclareMathOperator{\supp}{supp}

\newcommand{\set}[2]{ \left\{ #1 :\, #2 \right\} }
\newcommand{\seqq}[2]{ \left\langle #1 :\, #2\right\rangle }
\newcommand{\seq}[1]{ \left\langle #1 \right\rangle }

\title[Countable products of countable equivalence relations]{Strong ergodicity around
countable products of countable equivalence relations}
\author{Assaf Shani}
\date{\today}

\begin{document}

\maketitle
\begin{abstract}
This paper deals with countable products of countable Borel equivalence relations and equivalence relations ``just above'' those in the Borel reducibility hierarchy.
We show that if $E$ is strongly ergodic with respect to $\mu$ then $E^\mathbb{N}$ is strongly ergodic with respect to $\mu^\mathbb{N}$.
We answer questions of Clemens and Coskey regarding their recently defined $\Gamma$-jump operations, in particular showing that the $\mathbb{Z}^2$-jump of $E_\infty$ is strictly above the $\mathbb{Z}$-jump of $E_\infty$.
We study a notion of equivalence relations which can be classified by infinite sequences of ``definably countable sets''. 
In particular, we define an interesting example of such equivalence relation which is strictly above $E_\infty^\N$, strictly below $=^+$, and is incomparable with the $\Gamma$-jumps of countable equivalence relations.

We establish a characterization of strong ergodicity between Borel equivalence relations in terms of symmetric models, using results from \cite{Sha18}.
The proofs then rely on a fine analysis of the very weak choice principles ``every sequence of $E$-classes admits a choice sequence'', for various countable Borel equivalence relations $E$.
\end{abstract}

\section{Introduction}

Let $E$ be an equivalence relation on a Polish space $X$. Say that $E$ is \textbf{Borel} if it is a Borel subset of $X\times X$, and $E$ is \textbf{countable} if each $E$-class is countable.
Given a countable group $\Gamma$ and a Borel action $a\colon\Gamma\action X$, the corresponding \textbf{orbit equivalence relation $E_a$} is defined by $x E_a y$ if there is some $\gamma\in\Gamma$ sending $x$ to $y$.
The Feldman-Moore theorem states that any countable Borel equivalence relation is the orbit equivalence relation of some Borel action of a countable group.

Given Borel equivalence relations $E$ and $F$ on Polish spaces $X$ and $Y$ respectively, a Borel map $f\colon X\lto Y$ is a \textbf{homomorphism from $E$ to $F$}, denoted $f\colon E\lto_B F$, if for any $x,x'\in X$, if $x\Erel x'$ then $f(x)\Frel f(x')$.
A Borel map $f\colon X\lto Y$ is a \textbf{reduction of $E$ to $F$} if for all $x,x'\in X$, $x\Erel x' \iff f(x)\Frel f(x')$. 
A Borel homomorphism $f\colon E\lto_B F$ correspond to a ``Borel definable'' map between $E$-classes and $F$-classes. This map is injective if and only if $f$ is a reduction.
Say that $E$ is \textbf{Borel reducible to} $F$, denoted $E\leq_B F$, if there is a Borel reduction from $E$ to $F$.
This pre-order is used to measure the complexity of various equivalence relations.
Say that $E<_BF$ if $E\leq_B F$ and $F\not\leq_B E$.

Let $a$ be a Borel action of a countable group $\Gamma$ on a probability measure space $(X,\mu)$.
Recall that the action is said to be \textbf{ergodic (with respect to $\mu$)} if every $a$-invariant Borel set has either measure zero or measure one.
An equivalent condition is: any Borel $a$-invariant function from $X$ to $[0,1]$ is constant on a measure one set.
Similarly, the action is said to be \textbf{generically ergodic} if any $a$-invariant Borel set is either meager or comeager, equivalently, if any Borel $a$-invariant function from $X$ to $[0,1]$ is constant on a comeager set.
Note that a 
Borel $a$-invariant function is precisely a Borel homomorphism from $E_a$ to $=_{[0,1]}$, where $=_Y$ is the equality relation on $Y$.

\begin{defn}[See {\cite[Definition 6.5]{Kec18b}}]
Let $E$ and $F$ be Borel equivalence relations on Polish spaces $X$ and $Y$ respectively and let $\mu$ be a probability measure on $X$.
Say that $\mathbf{E}$ \textbf{is} $\mathbf{(\mu,F)}$-\textbf{ergodic} if for any Borel homomorphism $f\colon E\lto_B F$ there is a Borel $E$-invariant measure one set $A\subset X$ such that $f$ maps $A$ into a single $F$-class.
\end{defn}
We sometime say that $E$ is $F$-ergodic with respect to $\mu$.
Note that $E$ is ergodic with respect to $\mu$ if and only if $E$ is $(\mu,=_{[0,1]})$-ergodic.
The notion of $E_0$-ergodicity is also known as ``strong ergodicity'', as first defined by Jones and Schmidt \cite{JS87}.
Let $\mathbb{F}_2$ be the free group on two generators, $\nu$ the $\left(\frac{1}{2},\frac{1}{2}\right)$ measure on $\{0,1\}$.
Define $E_0$ on $\{0,1\}^\mathbb{Z}$ to be the orbit equivalence relation of the shift action $\Z\action \{0,1\}^\mathbb{Z}$, and $E_\infty$ as the orbit equivalent relation of the shift action $\F\action \{0,1\}^{\F}$.
Then $E_\infty$ is $(\nu^{\F},E_0)$-ergodic (see \cite{HK05}).
The notion of strong ergodicity is prevalent in the study of countable Borel equivalence relations. In order to show that $E$ is not Borel reducible to $F$, one often shows that $E$ is $F$-ergodic with respect to some measure (see \cite{AK00}).

When studying non countable Borel equivalence relations, Baire category arguments are often used rather than measure theoretic ones. In this case the notion of generic strong ergodicity is often useful.
Let $E$ and $F$ be Borel equivalence relations on Polish spaces $X$ and $Y$ respectively. 
Say that $\mathbf{E}$ \textbf{is generically} $\mathbf{F}$-\textbf{ergodic} if for any Borel homomorphism $f\colon E\lto_B F$ there is a Borel $E$-invariant comeager set $A\subset X$ such that $f$ maps $A$ into a single $F$-class.
Say that $E$ is \textbf{generically ergodic} if it is generically $=_{[0,1]}$-ergodic. If $E$ is induced by an action of some countable group, this is equivalent to the action being generically ergodic.


\subsection{Strong ergodicity for infinite products}\label{subsec;intro-ergodicity}

Given equivalence relations $E_n$ on Polish spaces $X_n$ respectively, let their \textbf{(full support) product} $\prod_n E_n$ be the equivalence relation on $\prod_n X_n$ defined by $x \mathrel{(\prod_n E_n)} y$ if $x(n) \mathrel{(E_n)} y(n)$ for all $n$. Let $E^\N$ be the product $\prod_{n} E$.
These equivalence relations arise naturally in the study of the Borel reducibility hierarchy, see for example \cite{HK97}.

Let the \textbf{finite support product} $\prod_n^{\textrm{fin}}E_n$ be the equivalence relation on $\prod_n X_n$ defined by $x\mathrel{(\prod_n^{\textrm{fin}}E_n)} y$ if $x\mathrel{(\prod_n E_n)} y$ and $x(n)=y(n)$ for all but finitely many $n$.
While the full support product of countable equivalence relation is no longer a countable equivalence relation, the finite support product is. The finite support product operation was studied by Kechris where he showed the following \cite[Lemma 4.2]{Kec18b}: suppose $E_n$ is a countable Borel equivalence relation and $\mu_n$ is a probability measure on $X_n$.
Then $\prod_n^{\textrm{fin}}E_n$ is ergodic with respect to $\prod_n \mu_n$ if and only if $E_n$ is ergodic with respect to $\mu_n$ for every $n$.
Note that if the finite support product is ergodic, then so is the full support product.
While strong ergodicity is not preserved under finite support products, we show that it is preserved under full support products
\begin{lem}[Corollary~\ref{cor;erg-iff-power-erg} below]\label{lem;product-strong-ergodicity}
Suppose $E_n$ is a countable Borel equivalence relation and $\mu_n$ is a probability measure on $X_n$. Let $F$ be a countable Borel equivalence relation. Then 
\begin{equation*}
    \prod_n E_n\textrm{ is }(\prod_n \mu_n,F)\textrm{-ergodic if and only if }E_n\textrm{ is }(\mu_n,F)\textrm{-ergodic for all }n.
\end{equation*}
\end{lem}
The non trivial direction is right to left, that is, showing that the infinite product is $F$-ergodic with respect to the product measure.
Let us note that full support products are necessary, even for $E_0$-ergodicity.
Let $X=\{0,1\}$, $\mu=\left(\frac{1}{2},\frac{1}{2}\right)$ and $E=X^2$ the equivalence relation with a single equivalence class.
Then $E$ is $(\mu,E_0)$-ergodic but $\prod_n^{\textrm{fin}}E$ is \textit{not} $(\mu^{\N},E_0)$-ergodic (as it is Borel bireducible with $E_0$).
More generally, if $\mu$ does not concentrate on a single element and $E$ is any equivalence relation, then $\prod_n^{\mathrm{fin}}E$ is not $(\mu^{\N},E_0)$-ergodic.
Thus a different approach is necessary for the lemma above, the proof of which appeals to set-theoretic definability in symmetric models.
The proof shows in general that for ergodic $E_n$, homomorphisms $\prod_n E_n\lto_B F$
are determined by homomorphisms defined on finite products $\prod_{n<m}E_n$, on a measure 1 set.

\subsection{The $\Gamma$-jumps of Clemens and Coskey}\label{subsec;intro-CCjumps}

Recently Clemens and Coskey \cite{CC19} defined new ``gentle'' jump operators. In particular, these yield new interesting equivalence relations ``just above'' the countable products of countable equivalence relations.
For a countable group $\Gamma$, Clemens and Coskey define the \textbf{$\Gamma$-jump of $E$}, $E^{[\Gamma]}$ on $X^\Gamma$, by $
    x \mathrel{E^{[\Gamma]}} y \iff (\exists \gamma\in\Gamma) (\forall \alpha\in\Gamma) x(\gamma^{-1}\alpha) \mathrel{E} y(\alpha).$
The $\Gamma$-jumps generalize the usual shift actions of countable groups. For example, the $\mathbb{Z}$-jump of $=_{\{0,1\}}$ is $E_0$ and the $\mathbb{F}_2$-jump of $=_{\{0,1\}}$ is $E_\infty$, where $=_{\{0,1\}}$ is the equality relation on $\{0,1\}$.


Clemens and Coskey show that for any countable group $\Gamma$, 
if $E$ is a generically ergodic countable Borel equivalence relation then $E^{[\Gamma]}$ is $E_\infty$-generically ergodic, so the $\Gamma$-jumps are genuinely different from countable products of countable equivalence relations.
Whether distinct groups $\Gamma$ produce non Borel bireducible jumps $E^{[\Gamma]}$ was left open. In particular, they ask if the equivalence relations $E_\infty^{[\F]}$ and $E_\infty^\N\times E_0^{[\Z]}$ are different than $E_\infty^{[\Z]}$.
In Section~\ref{sec;CCjumps} we show that $E_\infty^\N\times E_0^{[\Z]}$ and $E_\infty^{[\Z]}$ are $\leq_B$-incomparable (see Corollary~\ref{cor;Z-jump-Einfty-vs-product} and Proposition~\ref{prop;product-vs-Z-jump-Einfty}), and that $E_\infty^{[\Z]}$ is strictly below $E_\infty^{[\mathbb{F}_2]}$ (and even $E_\infty^{[\mathbb{Z}^2]}$) with respect to Borel reducibility (see Corollary~\ref{cor;distinct-Z^d-jumps}).
Moreover, we completely characterize strong ergodicity between $\Gamma$-jumps of countable Borel equivalence relations in terms of group theoretic properties.
\begin{thm}\label{thm;ergodicity-between-gamma-jumps}
Let $\Gamma$ and $\Delta$ be countable groups and $E$ a generically ergodic countable Borel equivalence relation.
The following are
equivalent:
\begin{enumerate}
    \item There is a subgroup $\tilde{\Delta}$ of $\Delta$, a normal subgroup $H$ of $\tilde{\Delta}$ and a group homomorphism from $\Gamma$ to $\tilde{\Delta}/H$ with finite kernel;
    \item $E^{[\Gamma]}$ is not generically $E_{\infty}^{[\Delta]}$-ergodic;
\end{enumerate}
\end{thm}
The same is also true for measures. That is, if $E$ is ergodic with respect to some measure $\mu$, then we may replace (2) with ``$E^{[\Gamma]}$ is not $(\mu^{\Gamma},E_\infty^{[\Delta]})$-ergodic''.
\begin{cor}\label{cor;distinct-Z^d-jumps}
Let $E$ be a generically ergodic countable Borel equivalence relation.
\begin{itemize}
    \item $E^{[\mathbb{Z}]}<_B E^{[\mathbb{Z}^2]}<_B E^{[\mathbb{Z}^3]}<_B...<_B E^{[\mathbb{Z}^{<\omega}]}<_B E^{[\mathbb{F}_2]}$;
    \item $E^{\Z}$ and $E^{[\mathbb{Z}_2^{<\mathbb{N}}]}$ are $\leq_B$-incomparable.
\end{itemize}
\end{cor}
\begin{remark}
This is in stark contrast to the situation with countable equivalence relations, in which case any action of a countable abelian group induces an equivalence relation which is Borel reducible to $E_0$ (see \cite{GJ15}), and any countable Borel equivalence relation is Borel reducible to $E_0$ on a comeager set (see \cite[Theorem 12.1]{KM04}).
A central tool in these and other arguments showing hyperfiniteness of countable Borel equivalence relations is the use of the ``Marker Lemma'' (see \cite[Lemma 6.7]{KM04}).
The main tool in the proof of Theorem~\ref{thm;ergodicity-between-gamma-jumps} is Lemma~\ref{lem;A-gamma-indiscernibles}, which shows a strong failure of the Marker Lemma in the context of a countable group acting on a set of $E_0$-classes (see Remark~\ref{remark;markers-lemma}).
This shows that some kind of appeal to the Marker Lemma is necessary to show the hyperfiniteness of, for example, $\Z^2$-actions. 
\end{remark}

\subsection{Classification by countable sequences of definably countable sets}
Recall that $=^+$ is the equivalence relation on $\R^\mathbb{N}$ defined so that the map $x\in\R^{\mathbb{N}}\mapsto \set{x(n)}{n\in\N}$ is a complete classification.
That is, $x\mathrel{=^+}y$ if and only if $ \set{x(n)}{n\in\N}$=$ \set{y(n)}{n\in\N}$.
The complete invariants are countable sets of reals.
However, given such arbitrary countable set of reals $A$, there is no enumeration of $A$ which is definable from $A$. 
On the other hand, suppose $E$ is a countable Borel equivalence relation, induced by a Borel action of a countable group $\Gamma$.
The map $x\mapsto\set{\gamma\cdot x}{\gamma\in\Gamma}$, sending $x$ to its orbit, is a complete classification.
In this case the invariants are countable sets of reals which can be enumerated in a simple way: given any $y$ in the orbit of $x$, $\seqq{\gamma\cdot y}{\gamma\in\Gamma}$ provides an enumeration of the orbit.

For a countable Borel equivalence relation $E$, $E^\N$ can be classified by sequences of countable sets of reals $A=\seqq{A_n}{n<\omega}$ such that for each $n$ there is an enumeration of $A_n$ definable from $A$ and elements in the transitive closure of $A$.
Next we consider equivalence relations which can be classified in such a way, by countable sequences of \textit{definably} countable sets of reals.

\begin{defn}\label{defn;pi-simple}
Let $
    D=\set{f\in(\R^\mathbb{N})^\mathbb{N}}{\forall n,i,j (f(n)(i)\textrm{ is computable from }f(n+1)(j))}.
$
Define the equivalence relation $\mathbf{E_\Pi}$ on $D$ by $x \mathrel{(E_\Pi)} y$ if for each $n$, $\set{x(n)(i)}{i\in\N}=\set{y(n)(i)}{i\in\N}$.
That is, $E_\Pi$ is $(=^+)^\mathbb{N}$ restricted to the domain $D$.
\end{defn}
Given $x\in D$ let $A_n^x=\set{x(n)(j)}{j\in\N}$.
The map $x\mapsto A^x=\seqq{A^x_n}{n\in\N}$ is a complete classification of $E_\Pi$.
For every $n$ and any $z\in A^x_{n+1}$, all the reals in $A^x_n$ are computable from $z$, so there is a definable enumeration of $A^x_n$ using $z$.
That is, given the sequence $A^x$ we can definably witness that each $A^x_n$ is a countable set.
We show that $E_\Pi$ is not a product of countable Borel equivalence relations, and is also different than the $\Gamma$-jumps.
\begin{thm}\label{thm;e-pi}
\begin{enumerate}
    \item $E_\infty^\N<_B E_\Pi<_B =^+$ and $E_\Pi$ is pinned;
    \item $E_\Pi\not\leq_B E_\infty^{[\Gamma]}$ and $E_0^{[\Gamma]}\not\leq_B E_\Pi$ for any infinite countable group $\Gamma$. 
\end{enumerate}
\end{thm}
See Definition~\ref{def;pinned} for the definition of pinned.
Part (1) is proved in Section~\ref{sec;pcp} and part (2) is proved in Section~\ref{sec;CCjumps}.
In Section~\ref{sec;pcp} we give a more general definition attempting to capture those equivalence relations which can be classified by countable sequences of definably countable sets of reals (Definition~\ref{def;pcp}), and show that $E_\Pi$ is maximal among those (Theorem~\ref{thm;maxpcp}).

\begin{remark}
The only previously known examples of equivalence relations between $E_\infty^\N$ and $=^+$ were the non-pinned equivalence relations constructed by Zapletal in \cite{Zap11}.
Clemens and Coskey \cite{CC19} note that $E_\infty^{[\Z]}$ is strictly above $E_\infty^\N$, strictly below $=^+$ and is pinned, thus is much closer to products of countable Borel equivalence relations.
By the results above the equivalence relation $E_\Pi$ also sits in this gap, and is incomparable with the $\Gamma$-jumps of countable Borel equivalence relations.
\end{remark}

\subsection{Very weak choice principles}
Using the developments in \cite{Sha18}, the results above are proved by first reformulating the questions in terms of symmetric models and choice principles. In particular we isolate an equivalent condition for strong ergodicity between equivalence relations (Lemma~\ref{lem;symm-model-ergod}).
The results then rely on studying the following weak choice principles.

\begin{defn}\label{def;cce}
Let $E$ be a countable equivalence relation on a Polish space $X$.
Then {\bf choice for countable sequences of $E$ classes}, abbreviated $\mathrm{CC}[E^\N]$, stands for the following statement:
Suppose $A=\seqq{A_n}{n\in \N}$ is a countable sequence of sets $A_n\subset X$ such that each $A_n$ is an $E$-class.
Then $\prod_n A_n$ is not empty.
That is, every $E^\N$-invariant admits a choice function.
\end{defn}

\begin{thm}[Theorem~\ref{thm;ergod-implies-sep} below]
Suppose $E$ and $F$ are countable Borel equivalence relations on Polish spaces $X$ and $Y$ respectively, and there is a Borel probability measure $\mu$ on $X$ such that $E$ is $(\mu,F)$-ergodic. 
Then there is a model in which $\mathrm{CC}[F^\N]$ holds yet $\mathrm{CC}[E^\N]$ fails.
\end{thm}
In particular, there are many pairs of countable Borel equivalence relations $E$ and $F$ such that $\mathrm{CC}[E^\N]$ and $\mathrm{CC}[F^\N]$ are independent.
A curious point here is that these models are constructed as intermediate extensions of a random real generic extension (using the measure $\mu^\N$).
Furthemore, these arguments will not work using a Cohen real, due to the fact that all countable Borel equivalence relations are hyperfinite on a comeager set.

Recall that {\bf countable choice for countable sets of reals}, abbreviated here as $\mathbf{CC[\R]^{\aleph_0}}$, states that any countable sequence $A=\seqq{A_n}{n\in\N}$ of countable sets of reals $A_n\subset \R$ admits a choice function.
This is a very weak choice principle, commonly studied in the literature (see \cite{choicecons}).
Over $\mathrm{ZF}$, for a countable equivalence relation $E$, $\mathrm{CC}[E^\N]$  follows from $\mathrm{CC[\R]^{\aleph_0}}$ 
(since any Polish space is Borel isomorphic to $\R$).

The proof that $E_\Pi$ is not Borel reducible to $E_\infty^\N$ and to the $\Gamma$-jumps relies on finding a model in which $\mathrm{CC}[\R]^{\aleph_0}$ fails yet for any countable Borel equivalence relation $E$, $\mathrm{CC}[E^\N]$ holds (Theorem~\ref{thm;Pisep}).

\subsection*{Acknowledgements}
The results in this paper are partially from my PhD thesis.
I would like to express my sincere gratitude to my advisor, Andrew Marks, for his guidance and encouragement, and for numerous informative discussions.

Furthermore, I wish to thank Clinton Conley, Yair Hayut, Alexander Kechris, Menachem Magidor and Jindrich Zapletal for helpful conversations.
I am also indebted to John Clemens and Sam Coskey for sharing with me an early version of their paper.

\section{Preliminaries}\label{sec;prelm}

\begin{fact}[Folklore]
Suppose $A$ is a set in some extension of $V$.
Then there exists a minimal transitive model of ZF containing $V$ and $A$, denoted $V(A)$.
\end{fact}

A familiar instance is when $A$ is a real, then more can be said:
\begin{fact}[Folklore]\label{fact;ext-by-real}
Suppose $A$ is a set of ordinals in some extension.
Then $V(A)$ satisfies the axiom of choice. 
\end{fact}

For the results in this paper it suffices to consider $V=L[r]$, for some real $r$. 
In this case $V(A)$ is the usual Hajnal relativized $L$-construction, $L(r,A)$.

Working in some extension of $V$, let $\mathrm{HOD}_{V,A}$ be the collection of all sets which are heredetarily definable using $A$, parameters from $V$, and parameters from the transitive closure of $A$.
Then $\mathrm{HOD}_{V,A}$ is a transitive model of ZF, extending $V$ and containing $A$.
In the examples considered in this paper it will not matter if one takes the minimal model $V(A)$ or the model $\mathrm{HOD}_{V,A}$.

Note that in $V(A)$, the model $\mathrm{HOD}_{V,A}$ must be everything.
We will use this below. 
That is, for any $X\in V(A)$, there is some formula $\psi$, parameters $\bar{a}$ from the transitive closure of $A$ and $v\in V$ such that $X$ is the unique set satisfying $\psi(X,A,\bar{a},v)$.
Equivalently, there is a formula such that $X=\set{x}{\varphi(x,A,\bar{a},v)}$.
In this case say that $\bar{a}$ is a support for $X$. We will be particularly interested in sets with empty support. That is, those definable from $A$ and parameters in $V$ alone.

\begin{lem}[Folklore]\label{lem;mutgen}
Suppose $x,y$ are mutually generic over $V$. Then $V[x]\cap V[y]=V$.
\end{lem}

The reader is referred to \cite{Gao09} or \cite{Kano08} for a discussion on equivalence relations which are classifiable by countable structures.
Consider for example $=^+$ on $\R^\N$, and the complete classification map $x\mapsto A_x=\set{x(i)}{i\in\N}$ (that is, $x=^+ y$ if and only if $A_x=A_y$).
The main point which is used below is that if $E$ is classifiable by countable structures, there is a complete classification map $x\mapsto A_x$ where $A_x$ is hereditarily countable, and the assignment is absolute.
That is, the statement $A_x=A$ cannot change in a forcing extension.

\begin{lem}[{see \cite[Lemma 3.7]{Sha18}}]\label{lem;correspondence}
Suppose $E$ and $F$ are Borel equivalence relations on $X$ and $Y$ respectively and $x\mapsto A_x$ and $y\mapsto B_y$ are classifications by countable structures of $E$ and $F$ respectively.
Let $V\subset M$ be ZF models and $A$ an $E$-invariant in $M$.
Assume further that $E$ is Borel reducible to $F$.
Then there is a $B$-invariant, definable from $A$ and parameters in $V$ alone, such that $V(A)=V(B)$.
\end{lem}
By ``$A$ is an $E$-invariant'' we only ask that there is some generic extension of $M$ in which $A=A_x$ for some $x$.
For example, if $E$ is $=^+$ then $A$ can be any set of reals.
\begin{proof}[Proof sketch]
Fix a Borel reduction $f$ of $E$ to $F$.
The set $B$ is defined as the unique set such that in any generic extension of $V(A)$, if there is an $x\in X$ with $A=A_x$, then $B=B_{f(x)}$.
Since $f$ is a reduction, $A$ can be defined from $B$ as the unique set such that for any $x$ in a generic extension, if $B=B_{f(x)}$ then $A=A_x$.
\end{proof}

Given an ideal $I$ of Borel subsets of $X$ let $P_I$ be the poset of all Borel $I$-positive sets, ordered by inclusion, $p$ extends $q$ if $p\subset q$. 
The reader is referred to \cite{Zap08} or \cite{ksz} for the definition and a discussion on proper ideals. For the results in this paper we only need to consider the meager ideal, in which case $P_I$ is Cohen forcing, or the null ideal in which case $P_I$ is Random forcing. 

The following lemma characterizes strong ergodicity between Borel equivalence relations, which are classifiable by countable structures, in terms of symmetric models.
The proof follows from \cite[Section 3]{Sha18}, where it is shown that a Borel homomorphism corresponds to a definable set in the relevant symmetric model.
\begin{lem}\label{lem;symm-model-ergod}
Suppose $E$ and $F$ are Borel equivalence relations on $X$ and $Y$ respectively and $x\mapsto A_x$ and $y\mapsto B_y$ are classifications by countable structures of $E$ and $F$ respectively.
Let $I$ be a proper ideal over $X$.
The following are equivalent.
\begin{enumerate}
    \item For every partial homomorphism $f\colon E\lto_B F$, defined on some $I$-positive set, $f$ maps an $I$-positive set into a single $F$-class;
    \item If $x\in X$ is $P_I$-generic over $V$ and $B$ is an $F$-invariant in $V(A_x)$ which is definable only from $A_x$ and parameters in $V$, then $B\in V$.
\end{enumerate}
\end{lem}
\begin{proof}[Proof sketch]
(1)$\implies$(2). If $x$ is $P_I$-generic and $B\in V(A_x)$ is definable using only $A_x$ and a parameter in $V$  then by \cite[Proposition 3.7]{Sha18} $B=B_{f(x)}$ for some partial Borel homomorphism $f\colon E\lto_B F$ defined on an $I$-positive set.
By (1) $f$ maps an $I$-positive set $C$ into a single $F$-class. Taking two mutually generic $x,y\in C$, we see that $B$ is in $V[x]\cap V[y]=V$.

(2)$\implies$(1). Given a homomorphism $f$, let $x$ be $P_I$-generic in the domain of $f$ and let $B=B_{f(x)}$.
By \cite[Lemma 3.6]{Sha18} $B\in V(A)$ is definable from $A$ and parameters in $V$ alone, and so $B\in V$ by (2).
Let $C$ be some $I$-positive set forcing that $B_{f(x)}=\dot{B}$. Then all $P_I$-generics in $C$ are mapped into a single $F$-class.

\end{proof}

We briefly recall the definition of pinned equivalence relations. See \cite[Definition 7.1.2]{Kano08} and \cite{Zap11}.
\begin{defn}\label{def;pinned}
Let $E$ be an analytic equivalence relation on a Polish space $X$. Let $P$ be a poset and $\tau$ a $P$-name. The pair $\left<P,\tau\right>$ is \textbf{a virtual $E$-class} if $P\times P$ forces that $\tau_{l}$ is $E$-equivalent to $\tau_{r}$, where $\tau_{l}$ and $\tau_{r}$ are the interpretation of $\tau$ using the left and right generics respectively.
A virtual $E$-class $\left<P,\tau\right>$ is \textbf{pinned} if there is some $x\in X$ from the ground model such that $P$ forces that $\tau$ is $E$-equivalent to $\check{x}$.
Finally, \textbf{$E$ is pinned} if any virtual $E$-class is pinned.
\end{defn}
Assume $E$ is a Borel equivalence relation and $x\mapsto A_x$ is a complete classification using hereditarily countable structures. 
Then a virtual $E$-class simply corresponds to a set $A$ in the ground model and a pair $\left<P,\tau\right>$ where $P$ forces that $A_{\tau}=A$.
To see this, let $G_l\times G_r$ be $P\times P$-generic and let $x_l,x_r$ be the interpretations of $\tau$ according to $G_l$,$G_r$ respectively. Then $A=A_{x_l}=A_{x_r}$ is in $V[G_l]\cap V[G_r]=V$.
So a virtual $E$-class is a set $A$ (in the ground model) which is forced to be the invariant of some real.
This virtual $E$-class is pinned if and only if $A$ is $A_x$ for some $x$ in the ground model.
Thus $E$ is pinned if and only if ``being an $E$-invariant'' is absolute for forcing extensions.
For example, the equivalence relation $=^+$ is not pinned, as the set of all reals $\R$ is not the image of any $x\in\R^{\mathbb{N}}$ in the ground model, but it is after collapsing the continuum to be countable.
On the other hand, suppose $E$ is a countable Borel equivalence relation induced by some action $a\colon\Gamma\curvearrowright X$ of a countable group $\Gamma$, and consider the classification $x\mapsto\Gamma\cdot x$.
Then ``being an orbit of the action $a$'' is absolute, and $E$ is pinned.
Similarly, if $E$ is countable then $E^{\mathbb{N}}$ is pinned: if $A=\seqq{A_n}{n<\omega}$ is a sequence of $E$-classes, then any $x\in\prod_nA_n$ satisfies that $A$ is the invariant of $x$.
\begin{remark}
While ZF proves that any countable Borel equivalence relation is pinned, and more generally that any $F_\sigma$ equivalence relation is pinned, $\mathrm{CC}[E^{\mathbb{N}}]$ is used to show that $E^\mathbb{N}$ is pinned. Furthermore this is necessary, as $E^{\mathbb{N}}$ is not pinned in the models considered in Section~\ref{sec;ctbl-power}.
Larson and Zapletal \cite{LZ19} also noticed the consistency of ZF with ``$E_0^\mathbb{N}$ is pinned''. They further study pinned equivalence relations in choiceless models, with a focus on models of DC.
\end{remark}

\subsection{Notation}
We use $\omega$ to denote the set of natural numbers $\N=0,1,2,...$.
For an equivalence relation $E$ on $X$ and a subset $A\subset X$, its $E$-saturation is defined by $[A]_E=\set{x\in X}{(\exists y\in A)xEy}$.
When no poset is involved, we say that a set $X$ is generic over $V$ if $X$ is in some forcing extension of $V$.
If $x$ is a real in some generic extension of $V$ then $x$ is in fact $P$-generic over $V$ for some poset $P$.
In this case we write $V[x]$ for $V(x)$.
For a formula $\psi$ and a model $M$, we denote the relativization of $\psi$ to $M$ by $\psi^M$.

Say that $E$ and $F$ are \textbf{Borel bireducible}, denoted $E\sim_B F$ if $E\leq_B F$ and $F\leq_B E$.
We write $=_X^+$ for the equivalence relation on $X^\omega$ identifying $x,y$ if they enumerate the same subset of $X$. For any Polish space $X$, $=_X^+\sim_B =^+$.

\section{Countable products of countable equivalence relations}\label{sec;ctbl-power}

In this section we consider countable powers of countable Borel equivalence relations. That is, equivalence relations of the form $E^\omega$ where $E$ is a countable Borel equivalence relation.
For notational simplicity we give a proof of Lemma~\ref{lem;product-strong-ergodicity} for powers only, the general proof is similar.
To each such equivalence relation $E^\omega$ we associated a choice principle $\cce$ (Definition~\ref{def;cce}), which states that any countable sequence of $E$-classes admits a choice function.
First we note that if $E$ is Borel reducible to $F$, then $\ccf$ implies $\cce$, over ZF.
More generally:

\begin{prop}
Let $E$ and $F$ be countable Borel equivalence relations on Polish spaces $X$ and $Y$ respectively.
If $E^\omega$ is Borel reducible to $F^\omega$ then $\ccf$ implies $\cce$, over ZF.
\end{prop}
\begin{proof}
Assume that $\ccf$ holds and fix a sequence $\seqq{A_n}{n<\omega}$ such that each $A_n$ is a $E$-class.
It remains to show that $\prod_n A_n\neq \es$.
By Lemma~\ref{lem;correspondence} there is an $F^\omega$-invariant $B$ such that $L(A)=L(B)$, $B=\seqq{B_n}{n<\omega}$ where each $B_n$ is an $F$-class.
Applying $\ccf$, there is some $y\in\prod_n B_n$.
Now $L(A)=L(B)\subset L[y]$ and the latter is a model of ZFC, so there is some $x\in\prod_n A_n$ in $L(y)$.
\end{proof}
Next we separate these choice principles.
\begin{thm}\label{thm;ergod-implies-sep}
Let $E$ and $F$ be countable Borel equivalence relations on $X$ and $Y$ respectively, $\mu$ a Borel probability measure on $X$ and suppose that $E$ is $(\mu,F)$-ergodic. Then there is a model in which $\ccf$ holds yet $\cce$ fails.
\end{thm}

\begin{cor}
 $\mathrm{CC}[E_\infty^\omega]$ is strictly stronger than $\mathrm{CC}[E_0^\omega]$.
\end{cor}

Adams and Kechris \cite{AK00} showed that there is a continuum size family $\mathcal{F}$ of countable Borel equivalence relations such that for any distinct $E,F\in\mathcal{F}$, $E$ is $F$-ergodic and $F$ is $E$-ergodic, with respect to some measures on their domains.
By Theorem~\ref{thm;ergod-implies-sep} we conclude:
\begin{cor}
There is a continuum size family $\mathcal{F}$ of countable Borel equivalence relations such that for any distinct $E,F\in\mathcal{F}$, $\cce$ and $\ccf$ are independent.
\end{cor}

Towards proving Theorem~\ref{thm;ergod-implies-sep} fix $E$, $F$ and $\mu$ as in the theorem.
By Feldman-Moore theorem we may fix a countable group $\Gamma$ and a Borel action $a \colon \Gamma \curvearrowright X$ such that $E=E_a$.
For a Borel probability measure $\nu$ on a Polish space $Y$, let $R(\nu)$ be the poset $P_{I(\nu)}$ where $I(\nu)$ is the ideal of $\nu$-measure zero sets. That is, all $\nu$-positive measure Borel subsets of $Y$ ordered by inclusion.
Let $x=\seqq{x_n}{n<\omega}\in X^\omega$ be a $R(\mu^\omega)$ generic over $V$, $A_n=[x_n]_E$ and $A=\seqq{A_n}{n<\omega}$ the corresponding $E^\omega$-invariant.

We work now in $V(A)$, recall the definitions from Section~\ref{sec;prelm}. Given $X\in V(A)$, there is some formula $\varphi$, a parameter $v\in V$ and finitely many parameters $\bar{a}$ from the transitive closure of $A$ (the support of $X$) such that $X=\set{z}{\varphi(z,A,\bar{a},v)}$. 
Since each $A_i$ (which is in the transitive closure of $A$) is definable from $A$, we may assume that $\bar{a}$ is contained in $\bigcup_i A_i$.
If $a\in A_i$ then $a$ is definable from $x(i)$.
Thus the support of $X$ can be taken to be of the form $\bar{a}=\seqq{x(i)}{i\in s}$ where $s\subset \omega$ is finite.

The following proposition establishes the basic symmetric-model analysis of $V(A)$ that will be used.
The proof follows a similar outline to that of an analogous property of the ``basic Cohen model'' (see \cite[Proposition 2.1]{Bla81}).
One difference is the required permutations, which are here the ones preserving $E^\omega$.
Furthermore, we are working with a Random real and not a Cohen real.
We note that the proposition holds for a Cohen real $x$ as well, with the proof slightly simpler.

\begin{prop}\label{prop;gen-blass}
Suppose $X\in V(A)$ and $X\subset V$, $s\subset \omega$ is finite and $\bar{x}=\seqq{x(i)}{i\in s}$ is a support for $X$.
Then $X\in V[\bar{x}]$.
\end{prop}
In particular, any real $b\in V(A)$ is in $V[x\rest n]$ for some $n<\omega$.
\begin{proof}
Let $\Gamma^{<\omega}$ be the group of all infinite sequences $\seqq{\gamma_i}{i<\omega}$ such that $\gamma_i\in\Gamma$ and $\gamma_i=1$ for all but finitely many $i$.
Fix $\varphi$ and $v\in V$ such that $X=\set{z}{\varphi(z,v,\bar{x},A)}$ in $V(A)$.
Let $\Omega=\omega\setminus s$. Given $y\in X^\Omega$ denote by $\bar{x}^\frown y$ the element of $X^\omega$ whose restriction to $\Omega$ is $y$ and its restriction to $s$ is $\bar{x}$.
For $p\in R(\mu^\omega)$, say that $p$ \textbf{agrees} with $\bar{x}$ if $\set{y\in X^\Omega}{\bar{x}^\frown y \in p}$ has positive $\mu^\Omega$-measure. That is, it is a condition in $R(\mu^\Omega)$.

Work now in the $R(\mu^\omega)$-generic extension $V[x]$.
If $z\in Z$ then $\varphi^{V(A)}(z,v,\bar{x},A)$ holds in $V[x]$, so there is a condition $p$ which agrees with $\bar{x}$ forcing $\varphi^{V(\dot{A})}(\check{z},\check{v},\dot{\bar{x}},\dot{A})$.
Similarly, if $z\notin Z$ there is some $p$ which agrees with $\bar{x}$ forcing $\neg\varphi^{V(\dot{A})}(\check{z},\check{x},\dot{\bar{x}},\dot{A})$.
We will show that $Z$ is defined in $V[\bar{x}]$ as the set of all $z$ such that there is some condition $p\in R(\mu^\omega)$ which agrees with $\bar{x}$ and forces that $\varphi^{V(\dot{A})}(\check{z},\check{v},\dot{\bar{x}},\dot{A})$. 
It suffices to show following: for any $z\in V$ there are no $p_0,p_1$ which agree with $\bar{x}$ such that $p_1\force \varphi^{V(\dot{A})}(\check{z},\check{v},\dot{\bar{x}},\dot{A})$ and $p_0\force \neg\varphi^{V(\dot{A})}(\check{z},\check{v},\dot{\bar{x}},\dot{A})$. 

For contradiction, assume we have $p_0,p_1$ as above.
Let $q_i=\set{y\in X^\Omega}{\bar{x}^\frown y\in p_i}$.
Fix a large enough countable model $M$ and let $\tilde{q}_i$ be the set of all $y\in q_i$ which are $R(\mu^\Omega)$-generic over $M[\bar{x}]$.
Note that $\tilde{q}_i$ has positive measure.
Since $E$ is $(\mu,F)$-ergodic, $E$ is in particular ergodic with respect to $\mu$, hence $\Gamma$ acts ergodically.
By \cite[Lemma 4.2]{Kec18a} the countable group $\Gamma^{<\omega}$ acts ergodically on $X^\omega$ (which we identify here with $X^\Omega$).
It follows that there is some  $g\in \Gamma^{<\omega}$ such that $(g^{-1}\cdot \tilde{q}_1)\cap \tilde{q}_0$ has positive measure.
In particular there is some $x'\in (g^{-1}\cdot \tilde{q}_1)\cap \tilde{q}_0$.
That is, both $x_0=\bar{x}^\frown x'$ and $x_1=\bar{x}^\frown g\cdot x'$ are $R(\mu^\omega)$-generics over $M$, and they agree on $\dot{A}$ and $\dot{\bar{x}}$.
Furthermore, $x_0$ extends $p_0$, thus $\varphi(z,v,\bar{x},A)$ fails in $M(A)$, but $x_1$ extends $p_1$, so $\varphi(z,v,\bar{x},A)$ holds in $M(A)$, a contradiction.
\end{proof}

\begin{cor}\label{cor;cce-fails}
In $V(\seqq{A_n}{n<\omega})$ there is no choice function for $\seqq{A_n}{n<\omega}$.
In particular, $\cce$ fails.
\end{cor}
\begin{proof}
Otherwise, there is a choice function $r\in\prod_n A_n$ which is in $V[x\rest n]$ for some $n$, by the lemma.
However, $r(n)$ is generic over $V[x\rest n]$, a contradiction.
\end{proof}
\begin{cor}\label{cor;countable-product-not-countable}
If $X\in V(A)$ is a set of subset of $V$, and $X$ is countable in $V(A)$, then $V(X)\neq V(A)$.
\end{cor}
\begin{proof}
Using an enumeration of $X$, $X$ can be coded as a subset of $V$ and therefore $X$ is in $V[x\restriction n]$ for some $n<\omega$.
\end{proof}
In particular if $X$ is an $E_\infty$-invariant in $V(A)$ (a countable set of reals), then $V(A)\neq V(X)$.
By Lemma~\ref{lem;correspondence} it follows that $E^\omega$ is not essentially countable.

We will see that the choice separation in Theorem~\ref{thm;ergod-implies-sep} corresponds to strong ergodicity between $E^\omega$ and $F^\omega$, rather than $E$ and $F$.
First we show that the first follows from the latter.

\begin{lem}\label{lem;inf-power-ergod}
Suppose $E$,$F$ are countable Borel equivalence relations on $X$ and $Y$ respectively. 
Let $\mu$ be a Borel probability measure on $X$ and assume that $E$ is $(\mu,F)$-ergodic.
Then $E^\omega$ is $(\mu^\omega,F)$-ergodic.
\end{lem}

\begin{cor}\label{cor;erg-iff-power-erg}
Suppose $E$,$F$, and $\mu$ are as above. 
Then $E$ is $(\mu,F)$-ergodic if and only if $E^\omega$ is $(\mu^\omega,F)$-ergodic if and only if $E^\omega$ is $(\mu^\omega,F^\omega)$-ergodic.
\end{cor}
\begin{proof}
Assume first that $E$ is $(\mu,F)$-ergodic, then $E^\omega$ is $(\mu^\omega,F)$-ergodic by the lemma above.
Given a homomorphism $f\colon E^\omega\lto F^\omega$ the projections $f_n(x)=f(x)(n)$ are homomorphisms from $E^\omega$ to $F$.
Thus each $f_n$ is constant on a $\mu^\omega$-conull set $A_n$, and so $f$ is constant on the conull set $A=\bigcap_n A_n$.

Conversely, assume that $E^\omega$ is $(\mu^\omega,F^\omega)$-ergodic and fix $f\colon E\lto F$. Then $f^\omega\colon X^\omega\lto Y^\omega$, defined by $f^\omega(x)(n)=f(x(n))$, is a homomorphism from $E^\omega$ to $F^\omega$.
By assumption, there is a conull set $A\subset X^\omega$ which $f^\omega$ sends to a single $F^\omega$-class, $\seqq{[y_i]_F}{i<\omega}$.
Let $A_0$ be the projection of $A$ to the first coordinate. Then $A_0$ is conull and for any $x\in A_0$, $f(x)\in [y_0]_F$.
\end{proof}
The proof of Lemma~\ref{lem;inf-power-ergod} will appeal to Proposition~\ref{prop;gen-blass} to reduce the problem to that of finite powers, in which case a direct measure theoretic argument works.
\begin{prop}\label{prop;fin-power-ergod}
Let $E$, $F$, $\mu$ be as in Lemma~\ref{lem;inf-power-ergod}.
Then $E^n$ is $(\mu^n,F)$-ergodic.
\end{prop}
\begin{proof}
Suppose that $E$ and $E'$ are countable Borel equivalence relations on $X$ and $X'$ respectively and are $F$-ergodic with respect to $\mu$ and $\mu'$ respectively. We show that $E\times E'$ is $(\mu\times\mu',F)$-ergodic. The proposition is then established inductively.
Fix a homomorphism $f\colon E\times E'\lto F$.
For $x\in X$ define $f_x\colon X'\lto Y$ by $f_x(x')=f(x,x')$.
For each $x$, $f_x$ is a homomorphism from $E'$ to $F$. By assumption there is a $\mu'$-measure 1 set $C_x\subset X'$ and $y_x\in Y$ such that $f_x(x') F y_x$ for any $x'\in C_x$.
If $x_1 E x_2$ then $y_{x_1} F y_{x_2}$, since for any $x'\in C_{x_1}\cap C_{x_2}$, $y_{x_1} F f(x_1,x') F f(x_2,x')F y_{x_2}$.

Let $D\subset X\times Y$ be the set of all pairs $(x,y)$ such that for any measure 1 set $C\subset X'$ there is some $x'\in C$ with $f_x(x') F y$.
$D$ is Borel and has countable $Y$-sections.
By Lusin-Novikov uniformization (see \cite{Kec95}) there is a Borel function $g\colon X\lto Y$ such that $(x,g(x))$ is in $D$ for all $x\in X$.
Note that $g(x) F y_x$ for any $x\in X$ and so $g$ is a homomorphism from $E$ to $F$.
Since $E$ is $(\mu,F)$-ergodic, there is a measure 1 set $C\subset X$ and $y\in Y$ such that $g(x) F y$ for all $x\in C$.
Let $A=\bigcup_{x\in C}\{x\}\times C_x$, a $\mu\times\mu'$-measure 1 subset of $X\times X$. We claim that $f(x,x') F y$ for any $(x,x')\in A$, which concludes the proof.
Indeed, $f(x,x')=f_x(x')\,F\,y_x\,F\,g(x)\,F\,y$.
\end{proof}

\begin{proof}[Proof of Lemma~\ref{lem;inf-power-ergod}]
Fix a homomorphism $f\colon E^\omega\lto F$.
Fix a large enough countable elementary submodel $M$ containing $f$. Let $x$ be $R(\mu^\omega)$-generic over $M$, $A_n=[x(n)]_{E}$ and $A=\seqq{A_n}{n<\omega}$.
By Lemma~\ref{lem;correspondence} the $F$-class of $f(x)$, $B=\set{y}{f(x)Fy}$ is in $M(A)$ and has empty support.
In particular, $f(x)\in M(A)$.
By Proposition~\ref{prop;gen-blass} $f(x)\in M[x\rest m]$ for some $m<\omega$, and so $B\in M[x\rest m]$.
Note that $M[x\rest m]=M(\seqq{A_i}{i<m})$.

Fix $\varphi$ and $v\in M$ such that $B=\set{y}{\varphi^{M(A)}(y,v,A)}$. We claim that $B$ is definable in $M(\seqq{A_i}{i<m})$ using only $\seqq{A_i}{i<m}$ and $v$.
This follows from Proposition~\ref{prop;gen-blass}. Viewing $M(\seqq{A_i}{i<m})$ as the ground model then $B$ is a subset of the ground model and is definable in $M(\seqq{A_i}{i<m})(\seqq{A_i}{i\geq m})$ using only $\seqq{A_i}{i\geq m}$ and the ground model parameters $v,\seqq{A_i}{i<m}$.

By \cite[Lemma 2.5]{Sha18} there is a homomorphism $g\colon E^m\lto F$ defined on a $\mu^m$-positive measure set $C$ containing $x\rest m$ such that $B=B_{g(x\rest m)}$.
By assumption and by Proposition~\ref{prop;fin-power-ergod} $g$ sends a full measure subset of $C$ to a single $F$-class.
W.l.o.g. assume the above statements are forced by the maximal condition.
Now the set $D$ of all $\mu^\omega$-generics $x\in X^\omega$ such that $x\rest m\in C$ has positive $\mu^\omega$-measure and for any $x\in D$, $f(x) F g(x\rest m)$ lies in a single $F$-class.
Finally, since $D$ is a set of generics its saturation $[D]_{E^\omega}$ is Borel (see \cite[Theorem 2.29]{ksz}), invariant and measure 1 (since $E^\omega$ is ergodic with respect to $\mu^\omega$).
\end{proof}

By Corollary~\ref{cor;cce-fails}, Lemma~\ref{lem;inf-power-ergod} and Corollary~\ref{cor;erg-iff-power-erg}, the following proposition will finish the proof of Theorem~\ref{thm;ergod-implies-sep}.

\begin{prop}
Let $E$, $F$ and $\mu$ be as above such that $E^\omega$ is $(\mu^\omega,F^\omega)$-ergodic.
Let $x\in X^\omega$ be $\mu^\omega$-Random real generic over $V$, $A_n=[x(n)]_E$ and $A=\seqq{A_n}{n<\omega}$.
Then $V(A)\models \ccf$.
\end{prop}
\begin{proof}
Suppose $B=\seqq{B_n}{n<\omega}\in V(A)$ is a sequence of $F$-classes. $B$ is definable using $A$, parameters in $V$ and some parameters $\bar{a}$ from $\bigcup_n A_n$.
If $\bar{a}=\emptyset$ then by Lemma~\ref{lem;symm-model-ergod} $B$ is in $V$, and therefore admits a choice function in $V$.
Generally, fix some $m$ such that $\bar{a}$ is in $V[x\rest m]$.
Work now with $V[x\rest m]$ as the base model, forcing $\seqq{A_n}{n\geq m}$ over it. A similar argument shows that $B\in V[x\rest m]$, which is a model of ZFC, and so $B$ admits a choice function.
\end{proof}
We finish this section with a simple remark about the choice principles $\cce$.
Fix a countable Borel equivalence relation $E$ on $X$ and let $a\colon\Gamma\curvearrowright X$ be an action of a countable group $\Gamma$ on $X$ generating $E$.
Note that if $\seqq{A_n}{n<\omega}$ is a sequence of $E$ classes, then any choice function $x\in\prod_n A_n$ codes a countable enumeration of $\bigcup_n A_n$ (using some fixed enumeration of $\Gamma$).
It follows that $\cce$ is equivalent to the formally stronger statement, that the union of countably many $E$-classes is countable. 
In particular, it follows that if $F\subset E$, then $\cce$ implies $\ccf$: 
given a sequence of $F$ classes $\seqq{B_n}{n<\omega}$, let $A_n=[B_n]_E$ be the corresponding $E$-class.
Now a well ordering of $\bigcup_n A_n$ gives a well ordering of $\bigcup_n B_n$.

\comment{
Whenever $\seq{A_n;\,n<\omega}$ is a sequence of subsets of $X$ such that each $A_n$ is contained in an orbit of $\Gamma$, then $\prod_n A_n\neq \es$.
\begin{proof}
Fix some enumeration $\Gamma=\seq{\gamma_i;\,i<\omega}$.
Let $B_n=[A_n]_\Gamma$. 
Applying $\mathrm{CC}[\R]^\Gamma$ to $\seq{B_n;\,n<\omega}$, there is $g\in\prod_n B_n$.
Define $f(n)$ to be $\gamma_i\cdot g(n)$ for the minimal $i$ such that $\gamma_i\cdot g(n)\in A_n$.
Then $f\in\prod_n A_n$.
\end{proof}
}


\section{Equivalence relations which can be classified by sequences of countable sets of reals}\label{sec;pcp}

The following definition attempts to capture those equivalence relations which can be classified by invariants of the form $\seqq{A_n}{n<\omega}$ where each $A_n$ is a subset of a Polish space and $A_n$ is definably enumerated using some elements in $\bigcup_n A_n$ as a parameter.

\begin{defn}\label{def;pcp}
Let $E$ be an equivalence relation such that the domain of $E$ is a subset of some product space $X=\prod_n X_n$.
$E$ is said to be a \textbf{Pointwise Countable Product (PCP) relation} if there are Borel equivalence relations $F_n$ on $X_n$ such that $E=\prod_n F_n\rest\dom E$, and for every $n$, for any $xEy$, $y(n)$ is $\Delta^1_1$ in $x(n+1)$.
\end{defn}
In this case, define $A^{E,x}_n\equiv\set{y(n)}{yEx}$, the projection to $X_n$ of the equivalence class of $x$ (which will be noted as $A^x_n$ when $E$ is unambiguous).
The map sending $x$ to $\seqq{A^x_n}{n<\omega}$ is a complete classification of $E$, using invariants which are countable sequences of countable subsets of a Polish space.
Furthermore, given an invariant $\seqq{A^x_n}{n<\omega}$, one can definably enumerate each $A^x_n$. That is, fix any element $z\in A^x_{n+1}$, then all elements in $A^x_n$ are $\Delta^1_1$ in $z$.

\begin{lem}\label{lem;weakpcp}
Suppose $E$ satisfies the following weakening of PCP:
There is a function $\varphi\colon\omega\lto\omega$ such that for all $n$ and $fEg$, $g(n)$ is $\Delta^1_1$ in $f(0),...,f(\varphi(n))$.
Then $E$ is Borel reducible to a PCP relation.
\end{lem}
\begin{proof}
Define $\psi\colon\omega\lto\omega$ inductively by $\psi(0)=0$ and $\psi(n+1)=\sup_{k\leq\psi(n)}\varphi(k)$.
Let $Y_n=\prod_{m\leq\psi(n)}X_n$ and consider the map $\theta\colon\prod_n X_n\lto \prod_n Y_n$ defined by $\theta(f)(n)=\seq{f(m);\,m\leq\psi(n)}$.
Let $G_n=\prod_{m\leq\psi(n)}F_n$ and $\tilde{E}=\prod_n G_n$.
It can be verified that $\theta$ is a reduction of $E$ to $\tilde{E}$, $\im (\theta)$ is Borel, and $\tilde{E}\rest\im(\theta)$ is a PCP relation.
\end{proof}

\begin{ex}
For any countable Borel equivalence relation $E$ on $X$, $E^\omega$ is Borel reduciblle to a PCP equivalence relation.
In this case, if $x,y\in X^\omega$ are $E^\omega$-related, $y(n)$ is $\Delta^1_1(x(n))$. 
\end{ex}


\begin{defn}
Given a PCP relation $E$ as above, let
\begin{equation*}
    A^E_n=\set{(y,z)}{\exists x\in \dom E(y=x(n+1)\wedge z\in A^{E,x}_n)}.
\end{equation*}
Note that in general it is analytic.
\end{defn}

\begin{lem}\label{lem;PCP-to-BorelPCP}
Suppose $E$ is a PCP Borel equivalence relation as in \ref{def;pcp} above.
So $E=\prod_n F_n\rest\dom E$ where $\dom E$ is a Borel subset of $\prod_n X_n$.
There is a Borel set $D$ containing $\dom E$ such that, for $F=\prod_n F_n\rest D$, $F$ is PCP relation and $A^F_n$ is Borel for all $n$.
\end{lem}
\begin{proof}
Given $x\in\prod_n X_n$ and $z\in X_n$, let $x[n,z]$ be the result of replacing the n'th coordinate of $x$ with $z$.
The idea is to determine whether $(x(n+1),z)\in A^E_n$ by asking if $x[n,z]$ is in $\dom E$.
For natural examples, this is in fact the case. In general, we will have to add members to $\dom E$.
To do this in a controlled manner, preserving the PCP conditions, we will use a reflection argument.

Define first $D_0$ to be all $x\in\prod_n X_n$ such that there is some $y\in\dom E$ for which $x\in\prod_n A^{E,y}_n$.
$D_0$ is $\Sigma^1_1$ and $\prod_n F_n\rest D_0$ is still PCP. 
That is, if $x(\prod_n F_n) y$ then $x(n)$ is $\Delta^1_1(y(n+1))$ for all $n$.
Define a property of subsets of $X$, $\Phi(A,B)$ as follows (we think of $A$ as $D_0$ and $B$ as its complement).
$\Phi(A,B)$ holds if
\begin{enumerate}
    \item For any $x,y$ if $x,y\in A$ and $x(\prod_n F_n)y$ then $x(n)$ is $\Delta^1_1(y(n+1))$ for all $n$;
    \item For any $x,y$, if $x\in A$ and for all $m$ there is some $x'\in A$ such that $x'(\prod_n F_n) x$ and $x'(m)=y(m)$, then $y\notin B$. 
\end{enumerate}
$\Phi$ is hereditary, continuous upward in the second variable and is $\Pi^1_1$ on $\Sigma^1_1$.
Furthermore, $\Phi(D_0,\neg D_0)$ holds.
By the second reflection theorem (Theorem 35.16 in \cite{Kec95}), there is some $\Delta^1_1$ set $D\supset D_0$ such that $\Phi(D,\neg D)$ holds.

Let $F=\prod_n F_n\rest D$. 
$F$ is PCP by condition (1) above.
By condition (2), for any $x\in D$ and $z\in A^{F,x}_n$, $x[n,z]\in D$.
Now the relation $(y,z)\in A^F_n$ holds if and only if $\exists x$ s.t. $x(n+1)=y$ and $x[n,z]\in D$, if and only if $\forall x$ if $x(n+1)=y$ then $x[n,z]\in D$.
Thus $A^F_n$ is $\Sigma^1_1$ and $\Pi^1_1$, and therefore is Borel. 

\end{proof}

\begin{prop}\label{prop;cover-borel-fucntions}
Suppose $E$ is PCP and Borel and $A^E_n$ are Borel.
There are Borel functions $\set{h^n_i}{i,n\in\omega}$ such that for any $x\in \dom E$, for all $n$, $A^{E,x}_n = \set{h^n_i(x(n+1))}{i\in\omega}$.
\end{prop}
\begin{proof}
For $x\in \dom E$, $\set{z}{(x(n+1),z)\in A^E_n}\subset A^{E,x}_n$ which is countable.
By Lusin-Novikov uniformization, there are Borel functions as desired.

\end{proof}

\begin{prop}\label{prop;pcp-is-pinned}
Assume $E$ is PCP and Borel, then 
\begin{enumerate}
\item $E\leq_B =^+$, and 
\item $E$ is pinned, hence is strictly below $=^+$.
\end{enumerate}
\end{prop}
\begin{proof}

By Lemma~\ref{lem;PCP-to-BorelPCP} and Proposition~\ref{prop;cover-borel-fucntions} we may assume the conclusion of Proposition~\ref{prop;cover-borel-fucntions}.
(1) Define a map $f\colon X\lto (\R^\omega)^\omega$ by $$f(x)=\seq{\seq{h_i^n(x(n+1));\,i<\omega};\,n<\omega},$$ 
sending $x$ to the sequence $\seqq{A^x_n}{n<\omega}$.
We show that $f$ is a Borel reduction to $(=^+)^\omega$ which is Borel bireducible with $=^+$.
If $xEy$ then $A^x_n=A^y_n$ for all $n$, and so $f(x)(n)=^+ f(y)(n)$ for each $n$.
Conversely, assume $f(x)(n)=^+f(y)(n)$ for every $n$.
For each $n$, there is some $i$ such that $x_n=h^n_i(x(n+1))$, and then some $j$ such that $h^n_i(x(n+1))=h^n_j(y(n+1))$, which is $F$-related to $y(n)$.
It follows that $x(n) F_n y(n)$ for each $n$, thus $xEy$.

(2) Note that the statement $\forall x\forall y(yEx\implies\exists i (y(n)=h^n_i(x(n+1)))$ is $\Pi^1_1$, and therefore is absolute.
Suppose $P$ is some poset such that $P\times P\force x_l E x_r$. Let $x$ be $P$-generic over $V$, we need to find some $z\in V$ such that $zEx$.
Take $y$ such that $(x,y)$ is $P\times P$ generic over $V$.

Using absoluteness between the models $V[x],\, V[y]$ and $V[x,y]$, and that $xEy$, it follows that $A^x_n=A^y_n$ for every $n$.
Thus $A^x=\seqq{A^x_n}{n<\omega}\in V[x]\cap V[y]=V$.
Applying countable choice in $V$, there is some $z$ such that $z(n)\in B_n$ for every $n$.
By absoluteness, there is some $x'\in V$ such that $x'\in\dom E$ and $x'(n) F_n z(n)$ for every $n$.
Thus $x'(n) F_n x(n)$ for each $n$, and $x'\in \dom E$, hence $x'Ex$ as required.
\end{proof}


\subsection{An interesting PCP equivalence relation}

We now define an interesting PCP relation which we denote $E_\Pi$. The definition below, which is the one used in all the proofs, is different than the one mentioned in the introduction (Definition~\ref{defn;pi-simple}). We will show in Proposition~\ref{prop;pi-simple} below that the two definitions are Borel bireducible.

\begin{defn}($E_\Pi$)\label{def;picomb}
Let $X_n=\R^{\omega^{n+1}}$, where we think of $X_{n+1}$ as $X_n^\omega$.
Let $X=\prod_n X_n$ and define:  \[
D=\set{f\in X}{\forall n \forall j(f(n+1)(j) \textrm{ is a permutation of }f(n))\wedge\exists j (f(n)=f(n+1)(j) )}
\]
We define the equivalence relation $\mathbf{E_\Pi}$ ($\Pi$ for product and permutations) on $D=\dom\ E_\Pi$ as follows:
$f E_\Pi g$ iff for every $n$, $\set{f(n)(i)}{i\in\omega}=\set{g(n)(i)}{i\in\omega}$.
That is, $E_\Pi$ is the equivalence relation $\prod_n (=_{\R^{\omega^n}}^+)$ restricted to the domain $D$.
\end{defn}
The equivalence relation $E_\Pi$ is PCP: if $f,g$ are $E_\Pi$-related, then for any $n$ there is some $j$ such that $g(n)=g(n+1)(j)$, and therefore there is some $i$ such that $g(n)=f(n+1)(i)$.
In particular, $g(n)$ is $\Delta^1_1(f(n+1))$.
The natural complete classification of $E_\Pi$ is the map sending $f\in\dom E_\Pi$ to the sequence $\seqq{A^{E_\Pi,f}_n}{n<\omega}$.
The complete invariants are sequences $\seqq{A_n}{n
<\omega}$ such that each $A_n$ is a subset of the Polish space $X_n$ and any element of $A_{n+1}$ is a countable enumeration of $A_n$.

By Proposition~\ref{prop;pcp-is-pinned}, $E_\Pi$ is pinned and strictly below $=^+$.
We now turn to prove that $E_\Pi$ is not Borel reducible to $E_\infty^\omega$.
The proof relies on constructing a model in which $\mathrm{CC}[\R]$ fails yet $\mathrm{CC}[E_\infty^\omega]$ holds.

Consider the poset $\set{p\colon\dom p\lto \omega}{\dom p\textrm{ is finite and }p\textrm{ is injective}}$.
This poset is isomorphic to Cohen forcing, and adds a generic permutation of $\omega$.

\begin{const}\label{const;Pi}
Let $P$ be the poset to add $\omega\times\omega$ mutually generic Cohen reals, indexed by $d_i$, $i<\omega$ where each $d_i$ is an $\omega$ sequence of Cohen reals.
We think of $d_0\in \R^\omega$ and $d_i$ for $i>0$ as a sequence of permutations of $\omega$.
Define inductively $a_n\in\R^{\omega^{n+1}}$  as follows:
\begin{itemize}
    \item Let $a_0=d_0\in \R^\omega$;
    \item Given $a_n$, let $a_{n+1}=\seq{a_n\circ d_{n+1}(i);\,i<\omega}$.
\end{itemize}
Thinking of $a_n$ as an element of $(\R^{\omega^{n}})^\omega$, let $A_n=\im a_n=\set{a_n(i)}{i\in\omega}$.
So $A_0\subset\R$ is a set of mutually generic Cohen reals.
$A_1\subset\R^\omega$ is a set of mutually generic enumerations of $A_0$. 
And so forth, $A_{n+1}$ is a set of generic enumerations of $A_n$.
Let $A=\seq{A_n;\,n<\omega}$. Our model will be $V(A)$.
\end{const}

\begin{prop}\label{prop;ccr0fail}
In $V(A)$, $\prod_n A_n=\es$.
In fact: 
Suppose $h\colon\omega\lto\omega$, in $V(A)$, is unbounded.
Then there is no function $f$ such that for each $n$, $f(n)$ is a non-empty finite subset of $A_{h(n)}$.
\end{prop}

\begin{thm}\label{thm;Pisep}
In V(A): For any countable Borel equivalence relation $E$, $\cce$ holds.
\end{thm}

It follows that $V(A)$ is not of the form $V(B)$ for any sequence of $E_\infty$-classes $B=\seqq{B_n}{n<\omega}$:
Otherwise, by the theorem there is some $x\in\prod_n B_n$ in $V(A)$.
As $B$ is definable from $x$, it follows that $V(A)=V(B)=V(x)$, which is a model of ZFC, by Fact~\ref{fact;ext-by-real}, contradicting the proposition above.
Note further that $A$ is an $E_\Pi$-invariant.
By Lemma~\ref{lem;correspondence} we conclude:

\begin{cor}\label{cor;e-pi-not-reducible}
$E_\Pi$ is not Borel reducible to $E_\infty^\omega$.
\end{cor}

The following lemma provides the basic symmetric model properties of $V(A)$ (the existence of minimal supports). We first show how to prove Theorem~\ref{thm;Pisep} from the lemma and then sketch a proof of the lemma.

\begin{lem}\label{lem;e-pi-model-basic}
There is map, definable in $V(A)$, sending a real $b\in V(A)$ to $n(b),E(b)$, where $n(b)<\omega$, $E(b)$ is the minimal finite subset of $A_n$ such that $b\in V(E(b))$ and $n$ is minimal for which such $E(b)$ exists.
\end{lem}

\begin{proof}[Proof of Proposition~\ref{prop;ccr0fail}]
Suppose $h,f$ are as in the statement of the proposition.
Using the linear ordering of the reals, define $g(n)$ to be the smallest member of $f(n)$.
$g$ is a real, so by Lemma~\ref{lem;e-pi-model-basic} there is some $m$ and a finite $E\subset A_m$ such that $g\in V(E)$.
In particular, $h(n)\leq m$ for all $n$.
\end{proof}

\begin{proof}[Proof of Theorem~\ref{thm;Pisep}]\label{pf;pisep}
Fix a countable Borel equivalence relation $F$ on a Polish space $X$.
Assume the parameters defining $F$ are in $V$ (otherwise, need to add a fixed finite set to all the supports below).
Assume $B=\seq{B_n;\,n<\omega}$ is a sequence of $F$-classes.
For $x,y\in B_n$, $y$ is $\Delta^1_1$ in $x$ (and a parameter for $E$). 
It follows that $x\in V(E)\iff y\in V(E)$ for any $E$, therefore $n(x)=n(y)$ and $E(x)=E(y)$.

Consider the map $h\colon \omega \lto \omega$ defined by 
$h(k)=n(x)$ for any $x\in A_k$,
and $f\colon\omega\lto\bigcup_n A_n$ defined by $f(k)=E(x)$ for any $x\in A_k$.
These are well defined by the argument above.
Note that, if $h(k)>0$, then $f(k)$ is a non-empty subset of $A_{h(k)}$.
By Proposition~\ref{prop;ccr0fail}, $h$ must be bounded.
Fix such bound $m\in \omega$,
then $B_n\subset V(A_m)$ for all $n$. 
Fix an enumeration $a$ of $A_m$ in $V(A)$. Then $V(A_m)\subset V[a]$, and $V[a]$ satisfies choice. In particular there is a well ordering in $V(A)$ of the reals in $V(A_m)$. Using this well order a choice function for $\seqq{B_n}{n<\omega}$ can be defined in $V(A)$.
\end{proof}


For $b\in V(A)$ and a finite $E\subset \bigcup_n A_n$, say that $E$ is a support for $b$ if there is a formula $\phi$ and $v\in V$ such that $\phi(b,v,E,A)$ defines $b$ in $V(A)$.

\begin{lem}\label{lem;E-Pi-supports}
For any real $b\in V(A)$, if $E$ is a support for $b$ then $b\in V(E)$.
\end{lem}
\begin{proof}
For notational simplicity consider the following special case:
assume that $E=a_1(0)\in A_1$, we show that $b\in V(a_1(0))$.
Let $\phi$ and $v\in V$ be such that $b=\set{n\in\omega}{V(A)\models \phi(n,v,a_1(0),A)}$.
First note that the sequence $\seqq{a_i}{1\leq i<\omega}$ can be added generically over $V(A_0)$.
Define $Q$ in $V(A_0)$ as the poset of all finite partial functions $p\colon\omega\times\omega\to A_0$ such that $p(k,\_)$ is injective.
Then $a_1$ is $Q$-generic over $V(A_0)$.
Let $P_2$ be the sub-poset of $P$ to add the sequence of permutations $\left<d_2,d_3,...\right>$. Then $\seq{a_1,d_2,d_3,...}$ is $Q\times P_2$-generic over $V(A_0)$ and $\seq{a_1,a_2,...}$ can be defined in $V(A_0)[\seq{a_1,d_2,d_3,...}]$.

We show that $b$ can be defined in $V(a_1(0))=V(A_0)[a_1(0)]$ as the set of all $n\in\omega$ such that any condition in $Q\times P_2$ which agrees with $a_1(0)$ forces $\phi^{V(\dot{A})}(n,v,\dot{a}_1(0),\dot{A})$.
The proof is similar to Proposition~\ref{prop;gen-blass} (see also \cite[Proposition 2.1]{Bla81} or \cite[Lemma 2.4]{Sha18}).
The main point is showing that two conditions which agree with $a_1(0)$ agree on $\phi^{V(\dot{A})}(n,v,\dot{a}_1(0),\dot{A})$.

Assume for contradiction that there are two conditions $p,q$ in $Q\times P_2$ which agree on $a_1(0)$ yet force incompatible statements about $\phi^{V(\dot{A})}(n,v,\dot{a}_1(0),\dot{A})$ for some $n$.
We may assume that $\seq{a_1,d_2,d_3,...}$ extends $p$. We will construct a generic $\seq{a'_1,d'_2,d'_3,...}$ which extends $q$, computes the same $\dot{A}$ and satisfies $a'_1(0)=a_1(0)$, which leads to a contradiction.

First, let $\pi_1$ be a finite permutation preserving $0$ such that $a'_1=a_1\circ\pi_1$ agrees with the restriction of $q$ to $Q$. $a'_1$ is $Q$-generic over $V(A_0)$.
Note that $a'_1\circ d_2(i)$ may no longer agree with $a_1\circ d_2(i)=a_2(i)$. That is, $\seq{a'_1,d_2,...}$ calculate the ``wrong'' $A_2$.
Let $d''_2(i)=\pi^{-1}\circ d_2(i)$, so $\seq{a'_1,d''_2,...}$ now calculates the correct $a_2$. Note that $d''_2$ is a generic sequence permutation over $V(A_0)$.
At this point $q$ may not agree with $d''_2$. Let $\pi_2$ be a finite permutation such that $d'_2=d''_2\circ\pi_2$ agrees with $q$. This is possible by genericity.
Note that $\seq{a'_1,d'_2,...}$ still calculates the correct $A_2$ (but not $a_2$).
Now set $d''_3(i)=\pi_2^{-1}\circ d_3(i)$, so that $\seq{a'_1,d'_2,d''_3,...}$ calculates the correct $a_3$. 
Continue is this fashion, at each step making finite changes to the values of $d_{k+1}(i)$ to get ``the correct $a_{k+1}$'' and then applying a finite permuting to the sequence $\seqq{d_{k+1}(i)}{i<\omega}$ to make it compatible with $q$.
Since $q$ has finite support, after finitely many steps we get $\seq{a'_1,d'_2,...,d'_k,d_{k+1},...}$ which is compatible with $q$, which completes the proof.
\end{proof}

If $E$ is a support for $b$, then for all large enough $n$, and any $a\in A_n$, all the elements of $E$ are definable from $a$, hence $\{a\}$ is a support for $b$.
Note that if $E$ is a support for $b$, $n$ is maximal such that $E\cap A_n\neq\es$, then $E\cap A_n$ is also a support for $b$.
Let $n(b)$ be the minimal $n$ such that there is some $E\subset A_n$ which is a support for $b$.

\begin{claim}
Fix a real $b\in V(A)$ and $n\in\omega$.
If $E_1,E_2\subset A_n$ are supports for $b$, then $E_1\cap E_2$ is a support for $b$.
By saying $\es\subset A_n$ is a support for $b$, when $n>0$, we mean that there is $E\subset A_{n-1}$ which is a support for $b$.
\end{claim}
\begin{proof}
The members of $A_n$ are enumerations of $A_{n-1}$, mutually generic over $V(A_{n-1})$.
By mutual genericity, if $b$ is in $V(A_{n-1})(E_1)$ and $V(A_{n-1})(E_2)$ then $b$ is in $V(A_{n-1})(E_1\cap E_2)$.
If $E_1\cap E_2=\es$, then $b\in V(A_{n-1})$ is definable from finitely many members of $A_{n-1}$, thus has a support in $A_{n-1}$.
\end{proof}

Let $E(b)$ be the minimal $E\subset A_{n(b)}$ which is a support for $b$. The map $b\mapsto n(b),\,E(b)$ satisfies the properties required by Lemma~\ref{lem;e-pi-model-basic}. This finishes the proof.

\subsection{More on $E_\Pi$}\label{subsec;maxpcp}
In this section we show that $E_\Pi$ is maximal Borel PCP.
\begin{thm}\label{thm;maxpcp}
If $E$ is a Borel PCP relation, then $E\leq_B E_\Pi$.
\end{thm}

\begin{proof}
Given a Polish space $Y$, let $E_\Pi(Y)$ be defined as in Definition~\ref{def;picomb} above, replacing each $\R^{\omega^n}$ with $Y^{\omega^n}$.
Any Borel isomorphism between $\R$ and $Y$ gives a Borel isomorphism between $E_\Pi=E_\Pi(\R)$ and $E_\Pi(Y)$.
Let $E$ be a $\mathrm{PCP}$ relation as in Definition~\ref{def;pcp} above, $\dom E\subset \prod_n X_n$.
The idea will be to construct a reduction of $E$ to $E_\Pi(X_0)$.
Each $f\in\dom E$ will be sent to an $E_\Pi(X_0)$-invariant $\seqq{D_n}{n<\omega}$ where $D_0$ is $A^f_0$, $D_1$ is a set of enumerations of $D_0$, given by the elements of $A^f_1$, and so forth.

By Lemma~\ref{lem;PCP-to-BorelPCP} and Proposition~\ref{prop;cover-borel-fucntions} we may fix Borel functions $h^n_i$ such that $A^f_n=\set{h_i^n(f(n+1))}{i\in\omega}$ for any $f\in \dom E$ and $n<\omega$.
Add two new and distinct elements, which we will call 0 and 1, to the space $X_0$.
Let $Y=X_0\cup\{0,1\}$, where $0,1$ are both isolated.
Fix Borel injections $\theta_n\colon X_{n+1}\lto \{0,1\}^\omega$ such that the constant sequences $\seq{0,0,...}$. and $\seq{1,1,...}$ are not in the images.

First we define some auxiliary functions:
given $x\in X_1$ define $\varphi_0(x)\in Y^\omega$ by:
\begin{itemize}
    \item $\varphi_0(x)(2i)=h_i^0(x)$;
    \item $\varphi_0(x)(2i+1)=\theta_0(x)(i)$.
\end{itemize}
Inductively, given $x\in X_{n+2}$ define $\varphi_{n+1}(x)\in ({Y^{\omega^{n+1}}})^\omega$ by
\begin{itemize}
    \item $\varphi_{n+1}(x)(2i)=\varphi_n(h^{n+1}_i(x))$;
    \item $\varphi_{n+1}(x)(2i+1)=\theta_{n+1}(x)(i)$.
\end{itemize}
(Here, the elements are in the space $Y^{\omega^{n}}$. By $0,1$ we mean the constant functions of such.)
Finally, define a map $\Psi\colon \dom E\lto \prod_n Y^{\omega^{n+1}}$ by
\begin{equation*}
    \Psi(f)(n)=\varphi_n(f(n+1)).
\end{equation*}
It remains to show that $\Psi$ is a reduction of $E$ to $E_\Pi(Y)$.
First we show that the range of $\Psi$ is included in the domain of $E_\Pi(Y)$.
For $f\in\prod_n X_n$ and $n\in\omega$,
we need to show that $\Psi(f)(n)=\Psi(f)(n+1)(j)$ for some $j$.
Since $E$ is PCP, $f(n+1)$ is $\Delta^1_1$ in $f(n+2)$, so there is some $j$ such that $h^{n+1}_j(f(n+2))=f(n+1)$.
Thus
\begin{equation*}
    \Psi(f)(n+1)(2j)=\varphi_{n+1}(f(n+2))(2j)=\varphi_n(h_j^{n+1}(f(n+2)))=\varphi_n(f(n+1))=\Psi(f)(n).
\end{equation*}

Next we show that $fEg\implies \Psi(f) E_\Pi\Psi(g)$.
For $n=0$, note that for any $fEg$: 
\begin{equation*}
    \im\Psi(f)(0)=\im \varphi_0(f(1))=A^f_0\cup\{0,1\}=A^g_0\cup\{0,1\}=\varphi_0(g(1))=\im\Psi(g)(0).
\end{equation*}
Given $fEg$, $n$ and $i$, we need to find some $j$ such that $\Psi(g)(n+1)(i)=\Psi(f)(n+1)(j)$.
For odd $i$, it follows by the choice of $\theta$.
By the PCP assumption, there is some $j$ such that $h^{n+1}_i(g(n+2))=h^{n+1}_j(f(n+2))$.
Thus
\begin{equation*}
\begin{split}
        &\Psi(g)(n+1)(2i)=\varphi_{n+1}(g(n+2))(2i)=\varphi_n(h_i^{n+1}(g(n+2))) \\
        &= \varphi_n(h_j^{n+1}(f(n+2)))=\varphi_{n+1}(f(n+1))(2j)=\Psi(f)(n+1)(2j).
\end{split}
\end{equation*}

Finally, it remains to show that if $\Psi(f)E_{\Pi} \Psi(g)$ then $fEg$.
Recall that, since $E$ is a PCP relation, it is a product of relations $F_n$ on $X_n$.
It suffices to show that $f(n)F_n g(n)$ for all $n$.
Assume $\Psi(f)E_\Pi\Psi(g)$ and fix $n\in\omega$.
Fix $i$ such that $\Psi(g)(n)=\Psi(g)(n+1)(i)$, and find a $j$ such that $\Psi(g)(n+1)(i)=\Psi(f)(n+1)(j)$.
It follows that both $i,j$ are even.
Let $j=2k$, then
\begin{equation*}
    \varphi_n(g(n+1))=\Psi(g)(n)=\Psi(f)(n+1)(j)=\varphi_{n+1}(f(n+2))(j)=\varphi_n(h^{n+1}_k(f(n+2)))
\end{equation*}
Note that $x$ is coded in the odd entries of $\varphi(x)$ using the functions $\theta$.
Thus from the equation above it follows that $g(n+1)=h_k^{n+1}(f(n+2))\in A^f_{n+1}$, hence $g(n+1)F_{n+1}f(n+1)$.
Furthermore, $g(0)\in\im\Psi(g)(0)=\im\Psi(f)(0)$, so there is some $j$ such that $g(0)=\Psi(f)(0)(j)=\varphi_0(f(1))(j)$.
It follows from the definition of $\varphi_0$ that $g(0)\in A^f_0$, thus $g(0)F_0 f(0)$.
We established that $g(n)F_n f(n)$ for all $n$, thus $gEf$ and the proof is done.
\end{proof}

We will use the theorem above to establish a few basic properties of the equivalence relation $E_\Pi$.

\begin{prop}\label{prop;pi-simple}
Let $D=\set{f\in(\R^\omega)^\omega}{\forall n,i,j (f(n)(i)\textrm{ is computable in }f(n+1)(j))}$.
Then $(=^+)^\omega \rest D\sim_B E_\Pi$.
That is, definitions \ref{def;picomb} and \ref{defn;pi-simple} agree.
\end{prop}

\begin{proof}
Fix recursive bijections $\phi_n\colon\R^{\omega^{n}}\lto\R$ where $\phi_0=\mathrm{id}$.
Define a map $\prod_n (\R^{\omega^{n}})^\omega \lto (\R^\omega)^\omega$ by $f\mapsto \seq{\phi_n\circ f(n);\,n<\omega}$.
This map is a reduction of $\prod_n (=_{\R^{\omega^n}})^+$ to $(=^+)^\omega$ and its image, when restricted to the domain of $E_\Pi$, is contained in $D$.
Thus $E_\Pi\leq_B (=^+)^\omega \rest D$.
Note that $(=^+)^\omega \rest D$ is Borel PCP, hence is Borel reducible to $E_\Pi$ by the Theorem~\ref{thm;maxpcp}.
\end{proof}

\begin{cor}
$(E_\Pi)^\omega\sim_B E_\Pi$.
\end{cor}
\begin{proof}
$(E_\Pi)^\omega$ can be represented as an equivalence relation with domain contained in a space $\prod_{m,n}X_{m,n}$, satisfying that:
if $f(E_\Pi)^\omega g$ then $g(m,n)$ is Borel in $f(m,n+1)$. 
Let $s\colon\omega\lto\omega\times\omega$ be the snake enumeration of $\omega\times\omega$, and let $X_n=X_{s(n)}$.
This gives an isomorphism of $\prod_n X_n$ to $\prod_{m,n}X_{m,n}$.
The pullback of $(E_\Pi)^\omega$ produces a relation which is Borel isomorphic to $(E_\Pi)^\omega$ and satisfies the conditions of Lemma~\ref{lem;weakpcp}.
\end{proof}

We conclude by noting the following generalization of PCP equivalence relations.
\begin{defn}
Let $E$ be an equivalence relation such that the domain of $E$ is a Borel subset of some product space $X=\prod_n X_n$.
$E$ is said to be PCP* if there are Borel equivalence relations $F_n$ on $X_n$ such that $E=\prod_n F_n\rest\dom E$, and for every $x$ and every $n$, the set $A^{E,x}_n\equiv\set{y(n)}{yEx}$ is countable.
\end{defn}
Any PCP equivalence relation is PCP*.
The proof of Proposition~\ref{prop;pcp-is-pinned} can be established for PCP* equivalence relations by similar arguments.
That is, any PCP* equivalence relation is pinned and strictly below $=^+$.
We do not know whether there is any PCP* equivalence relation which is not reducible to $E_\Pi$.

\section{Applications to the Clemens-Coskey jumps}\label{sec;CCjumps}

Recall the definition of the $\Gamma$-jumps from Section~\ref{subsec;intro-CCjumps}.
Let $E$ be a countable Borel equivalence relation on a Polish space $X$ and $\Gamma$ a countable group.
Given $x\in X^{\Gamma}$, for each $\gamma\in\Gamma$ let $A^x_\gamma=[x(\gamma)]_E$, the $E$-class of $x(\gamma)$.
For a $\Gamma$-indexed sequence $A'=\seqq{A'_\alpha}{\alpha\in\Gamma}$, define $\gamma\cdot A'=\seqq{A'_{\gamma^{-1}\alpha}}{\alpha\in\Gamma}$.
Define $$A^x=\set{\gamma\cdot\seqq{A_{\alpha}}{\alpha\in\Gamma}}{\gamma\in \Gamma}.$$
The map $x\mapsto A^x$ is a complete classification of $E^{[\Gamma]}$.

Fix a countable infinite group $\Gamma$.
Let $E$ be either ergodic with respect to a measure $\mu$ on $X$ or generically ergodic.
Let $x\in X^\Gamma$ be a $\mu^\Gamma$-Random generic, or Cohen generic, respectively. Consider $A=A^x$, its $E^{[\Gamma]}$-invariant.
We will study the model $V(A)$. Note that this model is equal to $V(\seqq{A_\gamma}{\gamma\in\Gamma})$. The latter is generated by a countable sequence of $E$-classes, which is the model studied in Section~\ref{sec;ctbl-power}.
By genericity the elements of $\set{A_\alpha}{\alpha\in\Gamma}$ are distinct.
In particular there is a well defined action of $\Gamma$ on $\set{A_\alpha}{\alpha\in\Gamma}$ defined by $\gamma\cdot A_\alpha = A_{\gamma^{-1}\alpha}$.
\begin{lemma}\label{lem;A-gamma-indiscernibles}
In $V(A)$, the elements $\set{A_\alpha}{\alpha\in\Gamma}$ are indiscernibles over $A$ and parameters in $V$.
\end{lemma}
\begin{remark}\label{remark;markers-lemma}
In our context the Marker Lemma (see \cite[Lemma 6.7]{KM04}) is manifested in the following way. 
Consider the shift action of $\Gamma$ on $[0,1]^\Gamma$ and let $x\in [0,1]^\Gamma$ be generic. The natural invariant will be the unordered set $\set{x(\gamma)}{\gamma\in\Gamma}$ together with the action of $\Gamma$, $\set{(\gamma,x(\alpha),x(\gamma^{-1}\alpha))}{\alpha,\gamma\in\Gamma}$.
Then the Marker Lemma provides arbitrarily sparse subsets of $\set{x(\gamma)}{\gamma\in\Gamma}$, definable using only the invariant and parameters in $V$.
On the other hand, when $\Gamma$ is acting on a collection of $E_0$-classes rather than on reals, Lemma~\ref{lem;A-gamma-indiscernibles} shows that the Marker Lemma fails in a strong way. 
\end{remark}
\begin{proof} 
We present the case where $x$ is Random generic. The case for a Cohen generic is similar and slightly easier.
Work in the big generic extension $V[x]$. Fix a formula $\phi$, $v\in V$ such that $\phi^{V(A)}(A_1,A,v)$ holds.
We will show that $\phi^{V(A)}(A_{\gamma^{-1}},A,v)$ holds for arbitrary $\gamma\in\Gamma$.
Assume towards a contradiction that $\phi^{V(A)}(A_{\gamma^{-1}},A,v)$ fails, and let $p$ be a condition forcing the above.
Since the shift action preserve the product measure, $\gamma\cdot p$ is a condition.
Furthermore, since $p$ forces $\phi^{V(\dot{A})}(\dot{A}_1,\dot{A},\check{v})$ then $\gamma\cdot p$ forces $\phi^{V(\dot{A})}(\dot{A}_{\gamma^{-1}},\dot{A},\check{v})$.

Consider $p$ and $\gamma\cdot p$ as positive measure subsets of $X^{\Gamma}$.
By \cite[Lemma 4.2]{Kec18a} the finite support power $\prod^{\mathrm{fin}}_{\gamma\in\Gamma} E$ is $\mu^{\Gamma}$-ergodic.
It follows that we may find generics $x_1,x_2$ such that $x_1$ extends $p$, $x_2$ extends $\gamma\cdot p$, and furthermore $x_1$ and $x_2$ are $E^\Gamma$-equivalent and they differ in only finitely many coordinates. (See the arguments in the proof of Proposition~\ref{prop;gen-blass}.)
Note that $A^{x_1}_\gamma=A^{x_2}_\gamma$ for every $\gamma$. Since $x_1$ extends $p$, working in $V[x_1]$ we conclude that $\phi^{V(A^{x_1})}(A^{x_1}_{\gamma^{-1}},A^{x_1},v)$ fails. However, since $x_2$ extends $\gamma\cdot p$, working in $V[x_2]$ we conclude that $\phi^{V(A^{x_2})}(A^{x_2}_{\gamma^{-1}},A^{x_2},v)$ holds, a contradiction.
\end{proof}

\begin{cor}\label{cor;no-A-gamma-definable}
No nonempty finite subset $\bar{A}\subset\set{A_\alpha}{\alpha\in\Gamma}$ is definable in $V(A)$ from $A$ and parameters in $V$ alone.
\end{cor}
\begin{proof}
Assume for contradiction there is such $\bar{A}$, defined as the unique solution to $\phi(\bar{A},A,v)$.
Fix $\gamma$ such that $\gamma\cdot\bar{A}\equiv\set{\gamma\cdot A_\alpha}{A_\alpha\in \bar{A}}$ is different than $\bar{A}$ (possible since $\Gamma$ is infinite.)
Then $\phi$ also uniquely defines $\gamma\cdot\bar{A}$, a contradiction.
\end{proof}

\begin{claim}\label{claim;Gamma-jump-erg}
If $B\in V(A)$ is a countable set of subsets of $V$, definable from $A$ and parameters in $V$ alone, then $B\in V$.
\end{claim}
\begin{proof}
Fix some enumeration $f$ of $B$. 
Note that $\supp(b)\subset \supp(f)$ for any $b\in B$.
Let $\bar{A}$ be the minimal finite subset of $\set{A_\alpha}{\alpha\in\Gamma}$ such that $\supp(b)\subset \bar{A}$ for any $b\in B$.
Since $\bar{A}$ is definable from $A$, $\bar{A}$ is empty by Corollary~\ref{cor;no-A-gamma-definable}.
It follows that $B$ is a subset of $V$. Then $B\in V$ by Proposition~\ref{prop;gen-blass}.
\end{proof}
By Lemma~\ref{lem;symm-model-ergod} we conclude:
\begin{cor}
Suppose $E$ is a countable Borel equivalence relation and either (1) $E$ is ergodic with respect to $\mu$ or (2) $E$ is generically ergodic. Then for any countable Borel equivalence relation $F$
\begin{enumerate}
    \item $E^{[\Gamma]}$ is $(\mu^\Gamma,F)$-ergodic;
    \item $E^{[\Gamma]}$ is generically $F$-ergodic.
\end{enumerate}
\end{cor}
Part (2) was proved by Clemens and Coskey \cite{CC19}. From part (1) we conclude:
\begin{cor}\label{cor;Z-jump-Einfty-vs-product}
$E_\infty^{[\Z]}$ is not Borel reducible to $E_\infty^\N\times E_0^{[\Z]}$.
\end{cor}
\begin{proof}
Fix a measure $\mu$ for which $E_\infty$ is $(\mu,E_0)$-ergodic.
Then $E_\infty^{[\Z]}$ is $(\mu^\Z,E_\infty)$-ergodic by the discussion above, and is therefore $(\mu^\Z,E_\infty^\N)$-ergodic.
Any reduction from $E_\infty^{[\Z]}$ to $E_\infty^\N\times E_0^{[\Z]}$ would give a reduction to $E_0^{[\Z]}$ on a $\mu^\Z$-measure 1 set.
Since $E_\infty$ is not reducible to $E_0^{[\Z]}$ on any $\mu$-measure 1 set (see \cite{CC19}), this would give a contradiction.
\end{proof}

\begin{prop}
$E_\Pi \not\leq E_\infty^{[\Gamma]}$ and $E_0^{[\Gamma]} \not\leq E_\Pi$ for any infinite countable group $\Gamma$.
\end{prop}
\begin{proof}
Let $A$ be an $E_0^{[\Gamma]}$-invariant of a Cohen generic real as above.
Suppose $B=\seq{B_0,B_1,...}$ is an $E_\Pi$-invariant definable from $A$ alone.
By Claim~\ref{claim;Gamma-jump-erg} $B_k\in V$ for every $k$, thus $B$ is in fact a subset of $V$, and so $B\in V$ by Proposition~\ref{prop;gen-blass}.
It follows that $E_0^{[\Gamma]}$ is generically $E_\Pi$-ergodic.

Let now $A$ be the $E_\Pi$-invariant from Construction~\ref{const;Pi}.
Let $B$ be a $E_\infty^{[\Gamma]}$-invariant in $V(A)$.
Then $V(B)=V(\seqq{B_\gamma}{\gamma\in \Gamma})$.
Since $V(A)$ is not generated by any countable sequence of $E_\infty$-classes (by Theorem~\ref{thm;Pisep}), it follows that $V(A)\neq V(B)$, so $E_\Pi\not\leq E_\infty^{[\Gamma]}$ by Lemma~\ref{lem;correspondence}.
\end{proof}

We now turn to prove Theorem~\ref{thm;ergodicity-between-gamma-jumps}, establishing strong ergodicity between the $\Gamma$-jumps for different values of $\Gamma$.
The central tools used in the proof are Proposition~\ref{prop;gen-blass} and Lemma~\ref{lem;A-gamma-indiscernibles}, which hold for both measure and category.
We focus on Baire category arguments, though the analogous results for measures also hold, as mentioned after Theorem~\ref{thm;ergodicity-between-gamma-jumps}. Note that when working mod meager sets, all the countable Borel equivalence relations below may be replaced with $E_0$, by generic hyperfiniteness \cite[Theorem 12.1]{KM04}.

\begin{lemma}
Let $\Gamma$ be a countable group and $\tilde{\Gamma}$ a finite normal subgroup of $\Gamma$. Let $E$ be countable Borel equivalence relation on $X$ which is generically ergodic.
Then $E^{[\Gamma]}$ is not generically $E_\infty^{[\Gamma/\tilde{\Gamma}]}$-ergodic.
\end{lemma}
\begin{proof}
Fix a Borel linear ordering $<$ of $X$ and let $m$ be the size of $\tilde{\Gamma}$.
Define an equivalence relation $\tilde{E}$ on $X^m$, relating two finite sequences $\bar{x},\bar{y}$ if there is a permutation $\pi$ of $m$ such that $x\circ\pi$ and $y$ are $E^m$-related.
Note that $\tilde{E}$ is a countable Borel equivalence relation, thus it suffices to show that $E^{[\Gamma]}$ is not generically $\tilde{E}^{[\Gamma/\tilde{\Gamma}]}$-ergodic.

Let $x\in X^\Gamma$ be Cohen generic and $A^x$ its $E^{[\Gamma]}$-invariant.
For $\gamma\in\Gamma$ define $B_{\gamma\tilde{\Gamma}}$ to be all choice functions in $\set{A_{\gamma\gamma'}}{\gamma'\in\tilde{\Gamma}}$. Note that $B_{\gamma\tilde{\Gamma}}$ can be viewed as an $\tilde{E}$-class.
Define $B=\set{\seqq{B_{g^{-1}\gamma\tilde{\Gamma}}}{\gamma\tilde{\Gamma}\in\Gamma/\tilde{\Gamma}}}{g\in\Gamma/\tilde{\Gamma}}$.
Then $B$ is an $\tilde{E}^{[\Gamma/\tilde{\Gamma}]}$-invariant definable in $V(A)$ from $A$ and parameters in $V$ alone, and $B\notin V$.
By Lemma~\ref{lem;symm-model-ergod} it follows that $E^{[\Gamma]}$ is not generically $\tilde{E}^{[\Gamma/\tilde{\Gamma}]}$-ergodic.

\end{proof}
\begin{cor}
Suppose there is a finite normal subgroup $\tilde{\Gamma}$ of $\Gamma$, a subgroup $\tilde{\Delta}$ of $\Delta$ and a normal subgroup $H$ of $\tilde{\Delta}$ such that $\Gamma/\tilde{\Gamma}$ embeds into $\tilde{\Delta}/H$.
Then $E^{[\Gamma]}$ is not generically $E_{\infty}^{[\Delta]}$-ergodic.
\end{cor}
\begin{proof}
Clemens and Coskey \cite{CC19} show that for any equivalence relation $F$, if $\Lambda$ is either a subgroup or a quotient of $\Gamma$, then $E^{[\Lambda]}$ is Borel reducible to $E^{[\Gamma]}$.
It follows that $E_\infty^{[\Gamma/\tilde{\Gamma}]}$ is Borel reducible to $E_\infty^{[\Delta]}$, so by the lemma above $E^{[\Gamma]}$ is not generically $E_\infty^{[\Delta]}$-ergodic.
\end{proof}
The corollary is the implication (1)$\implies$(2) of Theorem~\ref{thm;ergodicity-between-gamma-jumps}. The following proposition gives the reverse implication and will establish the theorem.
\begin{prop}
Let $\Gamma$ and $\Delta$ be countable groups, $E$ a generically ergodic countable Borel equivalence relation and $F$ a countable Borel equivalence relation.
Suppose $E^{[\Gamma]}$ is not generically $F^{[\Delta]}$-ergodic.
Then there is a subgroup $\tilde{\Delta}$ of $\Delta$, a normal subgroup $H$ of $\tilde{\Delta}$ and a group homomorphism from $\Gamma$ to $\tilde{\Delta}/H$ with finite kernel.
\end{prop}
\begin{proof}
By Lemma~\ref{lem;symm-model-ergod} there is a Cohen generic $x\in X^{\Gamma}$, $A=A^x$ its $E^{[\Gamma]}$-invariant, such that in $V(A)$ there is an $F^{[\Delta]}$-invariant $B$, definable from $A$ and parameters in $V$ alone, such that $B\notin V$. Let $B=\set{\seqq{B_{\delta^{-1}\xi}}{\xi\in\Delta}}{\delta\in\Delta}$, where each $B_\xi$ is an $F$-class.

\begin{claim}
There is some $\xi$ and $b\in B_\xi$ which is not in $V$.
\end{claim}
\begin{proof}
Otherwise, it follows that each $B_\xi$ is in $V$ and so $\seqq{B_{\delta^{-1}\xi}}{\xi\in\Delta}$ is a subset of $V$ for each $\delta$.
$B$ is countable in $V(A)$: fix some $B'=\seqq{B'_\xi}{\xi\in\Delta}\in B$. Let $\delta\cdot B'=\seqq{B'_{\delta^{-1}\xi}}{\xi\in\Delta}$, then $\seqq{\delta\cdot B'}{\delta\in\Delta}$ is an enumeration of $B$.
Thus $B$ is a countable set of subsets of $V$, definable from $A$ alone, so $B$ is in $V$ by Claim~\ref{claim;Gamma-jump-erg}.
\end{proof}
Fix some $b\in B_\xi$ not in $V$.
Since $b$ is a real there is a minimal finite $s\subset \Gamma$ such that for $\bar{x}=\seqq{x(\gamma)}{\gamma\in s}$, $b\in V[\bar{x}]$. 
Let $\bar{A}=\seqq{A_\gamma}{\gamma\in s}$.
Note that $V(\bar{A})=V[\bar{x}]$ and $B_\xi\subset V(\bar{A})$.
Say that $\bar{A}$ is the support for $B_\xi$.
For a fixed $\bar{A}$ there could be many $B_\xi$ whose support is $\bar{A}$.
We utilize the following coding functions to ensure that the set of $\xi$ for which $\bar{A}$ is a support for $B_\xi$ forms a subgroup. These are variations of the coding functions used by Clemens and Coskey in \cite{CC19} to show that $E^{[\mathbb{Z}]}$ is Borel reducible to $=^+$.

For $B'\in B$ define $p_{B'}\colon \Delta^2\times \Gamma\lto \{0,1\}$ by $p_{B'}(\delta_1,\delta_2,\gamma)=1$ if and only if the support of $(\delta_1\cdot B')_1$ is the $\gamma$-shift of the support of $(\delta_2\cdot B')_1$.
Let $P=\set{p_{B'}}{B'\in B}$. $P$ is a set of reals definable from $A$ alone (since $B$ is definable from $A$ alone).
Furthermore $P$ is countable:
fix any $p_{B'}\in P$ and define $p_\delta(\delta_1,\delta_2,\gamma)=1$ if and only if $p_{B'}(\delta_1\delta,\delta_2\delta,\gamma)=1$.
Given any other $p_{B''}\in B$, fix $\delta$ such that $B''=\delta\cdot B'$, then $p_{B''}(\delta_1,\delta_2,\gamma)=1$ if and only if the support of $(\delta_1\cdot B'')_1$ is a $\gamma$-shift of the support of $(\delta_2\cdot B'')_1$, if and only if the support of $(\delta_1\delta\cdot B')_1$ is a $\gamma$-shift of the support of $(\delta_1\delta\cdot B')_1$, if and only if $p_\delta(\delta_1,\delta_2,\gamma)=1$.
It follows that $\seqq{p_\delta}{\delta\in\Delta}$ enumerates $P$.
By Claim~\ref{claim;Gamma-jump-erg} we conclude that $P$ is in $V$, and so each $p\in P$ is in $V$.

Fix $p\in P$ and $B^\ast\in B$ such that $p_{B^\ast}=p$ and $\bar{A}$ is the support of $B^\ast_1$.
Let $\tilde{\Delta}$ be the set of all $\delta\in\Delta$ such that $p_{\delta\cdot B^\ast}=p$ and the support of $(\delta\cdot B^\ast)_1$ is $\gamma\cdot \bar{A}$ for some $\gamma\in\Gamma$.
Note that for any $\gamma$ there is a $\delta$ as above, by indiscernibility.
Define $H\subset\tilde{\Delta}$ as the set of all $\delta\in\tilde{\Delta}$ such that $B^\ast_1$ and $(\delta\cdot B^\ast)_1$ have the same support $\bar{A}$.
From now on $\bar{A}$, $B^\ast$ and $p$ are fixed. For $\delta\in\tilde{\Delta}$ say that \textbf{$\gamma$ is the support of $\delta$} if $\gamma\cdot\bar{A}$ is the support of $(\delta\cdot B^\ast)_1$
\begin{lemma}\label{lem;proof-of-group-homo}
Suppose $\delta,\delta'\in\tilde{\Delta}$ with supports $\gamma,\gamma'$ respectively.
Then $\delta^{-1}, \delta\delta'\in \tilde{\Delta}$ with supports $\gamma^{-1}, \gamma\gamma'$ respectively.
\end{lemma}
\begin{proof}
Suppose $\delta\in\tilde{\Delta}$ and $\gamma\cdot\bar{A}$ is the support for $(\delta\cdot B^\ast)_1$.
Then $p_{\delta\cdot B^\ast}(1,\delta^{-1},\gamma)=1$. Since $p_{B^\ast}=p=p_{\delta\cdot B^\ast}$, $p_{B^\ast}(1,\delta^{-1},\gamma)=1$ as well. Thus $\gamma^{-1}\cdot\bar{A}$ is the support for $(\delta^{-1}\cdot B^\ast)_1$.
Furthermore, $p_{\delta^{-1}\cdot B^\ast}(\delta_1,\delta_2,\gamma)=1$ if and only $p_{B^\ast}(\delta_1 \delta^{-1},\delta_2 \delta^{-1},\gamma)=1$ if and only if $p_{\delta\cdot B^\ast}(\delta_1 \delta^{-1},\delta_2 \delta^{-1},\gamma)=1$ if and only if $p_{B^\ast}(\delta_1 \delta^{-1}\delta,\delta_2 \delta^{-1}\delta,\gamma)=1$ if and only if $p_{B^\ast}(\delta_1,\delta_2,\gamma)=1$. Therefore $\delta^{-1}$ is in $\tilde{\Delta}$.

Assume now $\delta,\delta'$ are both in $\tilde{\Delta}$, where $\gamma\cdot\bar{A}$ and $\gamma'\cdot\bar{A}$ are the supports for $(\delta\cdot B^\ast)_1$ and $(\delta'\cdot B^\ast)_1$ respectively.
Then $p_{B^\ast}(\delta,1,\gamma)=1=p_{\delta'\cdot B^\ast}(\delta,1,\gamma)$, therefore the support of $(\delta\delta'\cdot B^\ast)_1$ is the $\gamma$-shift of the support of $(\delta'\cdot B^\ast)_1$, that is, it is $\gamma\gamma'\cdot\bar{A}$.
Furthermore,
\begin{equation*}
    p_{\delta\delta'\cdot B^\ast}(\delta_1,\delta_2,\zeta)=p_{\delta'\cdot B^\ast}(\delta_1\delta,\delta_2\delta,\zeta)=
    p_{B^\ast}(\delta_1\delta,\delta_2\delta,\zeta)=
    p_{\delta\cdot B^\ast}(\delta_1,\delta_2,\zeta)=
    p_{B^\ast}(\delta_1,\delta_2,\zeta).
\end{equation*}
Therefore $\delta\delta'$ is in $\tilde\Delta$.
\end{proof}

\begin{cor}
$\tilde{\Delta}$ is a subgroup of $\Delta$ and $H$ is a normal subgroup of $\tilde{\Delta}$.
\end{cor}
\begin{proof}
It follows from the lemma that $\tilde{\Delta}$ and $H$ are subgroups. For example if $h\in H$ then its support is $1\in\Gamma$. By the lemma $h^{-1}\in \tilde{\Delta}$ with support $1^{-1}=1$, that is, $h^{-1}\in H$.
We now show that $H$ is a normal subgroup of $\tilde{\Delta}$.
Suppose $h\in H$, $\delta\in\tilde{\Delta}$ with support $\gamma$.
By the lemma $h\delta\in\tilde{\Delta}$ with support $\gamma$. Applying the lemma again it follows that $\delta^{-1}h\delta\in\tilde{\Delta}$ with support $\gamma^{-1}\gamma=1$, thus $\delta^{-1}h\delta\in H$.
\end{proof}
For $\gamma\in\Gamma$ define $H_\gamma\subset\tilde{\Delta}$ as the set of all $\delta\in\tilde{\Delta}$ whose support is $\gamma$.
It follows from Lemma~\ref{lem;proof-of-group-homo} that each $H_\gamma$ is a coset of $H$ and that the map $\gamma\mapsto H_\gamma$  is a group homomorphism.
Furthermore, if $\gamma$ is in the kernel then $H_\gamma=H$ and so $\gamma\cdot\bar{A}=\bar{A}$. There could be only finitely many such $\gamma$'s, so the kernel is finite.
\end{proof}

We now show that $E_\infty^{[\Z]}$ and $E_0^{[\Z]}\times E_\infty^\omega$ are pairwise $\leq_B$-incomparable.
By Corollary~\ref{cor;Z-jump-Einfty-vs-product} it suffices to show the following irreducibility. Note that we do not have strong ergodicity in this case, with respect to either measure or category.

\begin{prop}\label{prop;product-vs-Z-jump-Einfty}
$E_0^{[\Z]}\times E_0^\omega$ is not Borel reducible to $E_\infty^{[\Z]}$.
\end{prop}
\begin{proof}
Since we are dealing with an additive group we will write the action in an additive way: $A_\alpha+\gamma=A_{\alpha+\gamma}$.
Towards a contradiction, assume that there is a Borel reduction of $E_0^{[\Z]}\times E_0^\omega$ to $E_\infty^{[\Z]}$.
Let $(x,y)\in (2^\omega)^\Z\times(2^{\omega})^\omega$ be Cohen generic. 
Let $A,B$ be the $E_0^{[\Z]}\times E_0^\omega$-invariant.
That is, $A=A^x$ as above and $B=\seqq{[y(n)]_{E_0}}{n<\omega}$.
By Lemma~\ref{lem;correspondence} there is an $E_\infty^{[\Z]}$-invariant $C\in V(A,B)$ such that $V(A,B)=V(C)$ and $C$ is definable from $A,B$ and parameters in $V$.
Let $C_k$, $k\in \Z$ be $E_\infty$-classes such that $C=\set{\seqq{C_{n+k}}{n\in\Z}}{k\in\Z}$. 
\begin{claim}
\begin{enumerate}
    \item The members of $\set{A_n}{n\in\Z}$ are indiscernible over $A$ and parameters from $V(B)$.
    \item For any $X\subset V$ there is a minimal support $(s,t)$ where $s$ is a finite subset of $\Z$, $t$ a finite subset of $\omega$ such that $X\in V(\seqq{x(k)}{k\in s},\seqq{y(n)}{n\in t})=V(\seqq{A_k}{k\in s},\seqq{B_n}{n\in t})$.
\end{enumerate}
\end{claim}
Part (1) follows from Lemma~\ref{lem;A-gamma-indiscernibles}, working over $V(B)$ as the ground model.
Part (2) is proved similarly, working over $V(B)$ we get that $X$ is in $V(B)[\bar{x}]$. Now working over $V[\bar{x}]$ we get that $X$ is in $V[\bar{x}][\bar{y}]$.
In particular, any $c\in C_k$ has a minimal support $(s,t)$. The support is the same for any other $c'\in C_k$, and so we say that $(s,t)$ is the support of $C_k$ as well.
Given such $s,t$ we will denote $\bar{x}=\seqq{x(k)}{k\in s}$, $\bar{A}=\seqq{A_k}{k\in s}$, $\bar{y}=\seqq{y(n)}{n\in t}$ and $\bar{B}=\seqq{B_n}{n\in t}$.
We sometimes call the pair $(\bar{x},\bar{y})$ or $(\bar{A},\bar{B})$ the support of $C_k$ as well.

For each $C'\in C$ ($C'=\seqq{C'_k}{k\in\Z}$) define $p_{C'}\colon \Z^2\lto \{0,1\}$ by $p_{C'}(t,l)=1$ if and only if $C'_t$ and $C'_l$ have the same support.
Let $P=\set{p_{C'}}{C'\in C}$.
Then $P$ is a countable set of reals, definable from $(A,B)$ alone, and therefore $P\in V(B)$.
It follows that $P\in V(\seqq{B_n}{n\leq m})$ for some $m$.
By replacing $V$ with $V(\seqq{B_n}{n\leq m})$ and $B$ with $\seqq{B_n}{n\geq m}$, we may assume that $P\in V$.

Fix $p\in P$ in $V$ for which there is some $C'\in C$ with $p=p_{C'}$ and the support of $C'_0$ is $(\bar{A},\bar{B})$.
Suppose there is some $C''\in C$ such that $p_{C''}=p$ and $C''_0$ has the same support $(\bar{A},\bar{B})$.
Fix $m$ such that $C''= C'+m$, then for every $k\in\Z$, $C'_{mk}$ has the same support $(\bar{A},\bar{B})$.
That is, we associate to $(\bar{A},\bar{B})$ an arithmetic sequence in the $\Z$-ordering on $\set{C_n}{n\in\Z}$.
If not such $C''$ exists then this arithmetic sequence is a singleton ($m=0$).
Note that the for distinct pairs $(\bar{A},\bar{B})$ and $(\bar{A}',\bar{B}')$ the corresponding arithmetic sequences must be disjoint.
\begin{claim}
If $\bar{A}$ is not empty then for any $\bar{B}$ the arithmetic sequence corresponding to $(\bar{A},\bar{B})$ is a singleton (that is, $m=0$).
\end{claim}
\begin{proof}
Otherwise, there is an arithmetic sequence with common difference $m$ corresponding to $(\bar{A},\bar{B})$ with $m>0$.
By indiscernibility, for any $t$ there is an arithmetic sequence with common difference $m$ corresponding to $(\bar{A}+t,\bar{B})$.
Since there could be only finitely many disjoint arithmetic sequences with fixed common difference $m$, we arrive at a contradiction.
\end{proof}

\textbf{Case 1} Suppose there is a support $(\bar{A},\bar{B})$ for some $C_t$ where $\bar{A}$ is not empty. 
By the claim $(\bar{A},\bar{B})$ defines uniquely some $C'\in C$ with $\supp C'_0=(\bar{A},\bar{B})$. As before, there is some $k\in\Z$ such that $\set{(\bar{A}+l,\bar{B})}{l\in\Z}$ corresponds to the $k$-arithmetic sequence $\set{C'_{kl}}{l\in\Z}$ with $\supp C'_{kl}=(\bar{A}+l,\bar{B})$.
\begin{claim}
\begin{enumerate}
    \item Assume that $(\es,\bar{B}')$ is a support corresponding to an $m$-arithmetic sequence. Then $0<m\leq k$.
    \item Assume $\bar{A}'$ is not empty any $(\bar{A}',\bar{B}')$ is a support corresponding to a singleton, so $\set{(\bar{A}'+l,\bar{B}')}{k\in\mathbb{Z}}$ corresponds to an $m$-arithmetic sequence for some $m$. 
    Then $m\leq k$.
\end{enumerate}
\end{claim}
\begin{proof}
We prove (1), the proof of (2) is similar.
Fix $l$ such that for the unique $C'$ corresponding to $(\bar{A}+l,\bar{B})$, there is $0<j<k$ with $\supp C'_j=(\es,\bar{B}')$.
By indiscernibility, this is true for $(\bar{A}+l,\bar{B})$ for any $l$. Fixing $C'$, it follows that for any $l$ there is some $0<j<k$ with $\supp C'_{kl+j}=(\es,\bar{B}')$.
We conclude that $m>0$ and $m\leq k$.
\end{proof}
There could be at most finitely many disjoint $m$-arithmetic sequences with $0<m\leq k$. 
Therefore there are only finitely many $\bar{B}'$ for which there is some $\bar{A}'$ (possibly empty) so that $(\bar{A}',\bar{B}')$ is a support of some $C_t$.
It follows that there is some $m<\omega$ such that $C_t\in V(A,\seqq{B_n}{n<m})$ for all $t$.
Forcing now $\seqq{B_n}{n\geq m}$ over $V(A,\seqq{B_n}{n<m})$, we get that $C$ is a set of subsets of the ground model $V(A,\seqq{B_n}{n<m})$, which is countable and definable from $\seqq{B_n}{n\geq m}$ and parameters in $V(A,\seqq{B_n}{n<m})$. 
By Corollary~\ref{cor;countable-product-not-countable} it follows that $V(A,\seqq{B_n}{n<m})(\seqq{B_n}{n\geq k})\neq V(A,\seqq{B_n}{n<m})(C)$, contradicting the assumption that $V(C)=V(A,B)$.

\textbf{Case 2} For any $C_k$ its support is of the form $(\es,\bar{B})$ for some $\bar{B}$.
It follows that $C\in V(B)$, a contradiction.
\end{proof}
The proposition also implies that $E^{[\Z]}<_B (E^{[\Z]})^2$ for any generically ergodic countable Borel equivalence relation $E$.
Similar arguments show that $E^{[\Z]}<_B (E^{[\Z]})^2<_B(E^{[\Z]})^3<_B...$.


\begin{thebibliography}{1}
    \bibitem[AK00]{AK00} S. Adams and A.S. Kechris, Linear algebraic groups and countable Borel equivalence relations, J. Amer. Math. Soc., 13(4) (2000), 909-943.

    \bibitem[Bla81]{Bla81} Blass, Andreas. The Model of Set Theory Generated by Countably Many Generic Reals. The Journal of Symbolic Logic, vol. 46, no. 4, 1981, pp. 732-752.
    
    \bibitem[CC$\infty$]{CC19} John D. Clemens and Samuel Coskey, New Jump Operators on Borel Equivalence Relations. In preparation.

    \bibitem[DJK94]{DJK94} R. Dougherty, S. Jackson and A.S. Kechris, The structure of hyperfinite Borel equivalence relations, Trans. Amer. Math. Soc., 341(1) (1994), 193-225.
    
    
    
    
    
    \bibitem[Gao09]{Gao09} S. Gao, Invariant descriptive set theory, CRC Press, 2009. 

    \bibitem[GJ15]{GJ15} Su Gao and Steve Jackson. Countable abelian group actions and hyperfinite equivalence relations. Invent. Math., 201(1):309-383, 2015.    
    
    
    
    \bibitem[HK97]{HK97} Hjorth, Greg; Kechris, Alexander S. New Dichotomies for Borel Equivalence Relations. Bull. Symbolic Logic 3 (1997), no. 3, 329--346.

    \bibitem[HK05]{HK05} G. Hjorth and A. S. Kechris, Rigidity theorems for actions of product groups and countable Borel equivalence relations. Mem. Amer. Math. Soc. 177 (2005), no. 833.    
    
    
    \bibitem[HR98]{choicecons} P. Howard and J.E. Rubin, Consequences of the axiom of choice, Mathematical Surveys and Monographs, no. v. 1; v. 59, American Mathematical Society, 1998.

    \bibitem[Jec73]{jechchoice} Thomas J. Jech, The axiom of choice, North-Holland Publishing Co., Amsterdam, 1973, Studies in Logic and the Foundations of Mathematics, Vol. 75.

    \bibitem[JS87]{JS87} V.F.R. Jones and K. Schmidt, Asymptotically invariant sequences and approximate finiteness, Amer. J. Math., 109 (1987), 91-114.
    
    
    \bibitem[Kano08]{Kano08} V. Kanovei, Borel equivalence relations, Amer. Math. Soc., 2008.    
    
    \bibitem[KSZ13]{ksz} Kanovei V., Sabok M., Zapletal J.: Canonical Ramsey theory on Polish Spaces. Cambridge University Press, Cambridge (2013).

    
    \bibitem[Kec95]{Kec95} A. S. Kechris, Classical Descriptive Set Theory, Springer, 1995.
    
    \bibitem[Kec$\infty$a]{Kec18a} A. S. Kechris, Quasi-invariant measures for continuous group actions [to appear]
    
    \bibitem[Kec$\infty$b]{Kec18b} A. S. Kechris, The theory of countable Borel equivalence relations, preprint.
    
    
    \bibitem[KM04]{KM04} A.S. Kechris and B.D. Miller, Topics in orbit equivalence, Springer, 2004.
    
    \bibitem[LZ$\infty$]{LZ19} Paul Larson and Jindrich Zapletal, Geometric set theory. In preparation.
    
    \bibitem[Mos39]{Mos39} Mostowski, Andrzej. "Über die Unabhängigkeit des Wohlordnungssatzes vom Ordnungsprinzip." Fundamenta Mathematicae 32.1 (1939): 201-252. 
    
    \bibitem[Sha$\infty$]{Sha18} Assaf Shani, Borel reducibility and symmetric models. preprint     
    
    \bibitem[Zap08]{Zap08} Jindrich Zapletal. Forcing Idealized. Cambridge Tracts in Mathematics 174. Cambridge University Press, Cambridge, 2008.
    
    \bibitem[Zap11]{Zap11} Jindrich Zapletal. Pinned equivalence relations. Mathematical Research Letters, 18:1149-1156, 2011.
    
\end{thebibliography}
\end{document}